\documentclass[11pt,letterpaper]{amsart}

\usepackage[letterpaper,margin=.95in]{geometry}
\usepackage[all]{xy} 
\usepackage{amsmath, amssymb, amsfonts, latexsym, mdwlist, amsthm, amscd,mathabx,braket}
\usepackage{subfig}
\usepackage{graphicx}
\usepackage{wrapfig}
\usepackage{etex}
\usepackage{tikz-cd}
\usepackage[colorlinks=true, citecolor=blue, urlcolor=blue, linkcolor=blue, pagebackref]{hyperref}

\usepackage{tikz}
\usetikzlibrary{calc,trees,positioning,arrows,chains,shapes.geometric,%
    decorations.pathreplacing,decorations.pathmorphing,shapes,%
    matrix,shapes.symbols}

\tikzset{
>=stealth',
  punktchain/.style={
    rectangle,
    rounded corners,
    draw=black, thick,
    minimum height=3em,
    text centered,
    on chain},
  line/.style={draw, thick, <-},
  eLement/.style={
    tape,
    top color=white,
    bottom color=blue!50!black!60!,
    minimum width=8em,
    draw=blue!40!black!90, very thick,
    text width=10em,
    minimum height=3.5em,
    text centered,
    on chain},
  every join/.style={->, thick,shorten >=1pt},
  decoration={brace},
  tuborg/.style={decorate},
  tubnode/.style={midway, right=2pt},
}

\usepackage{paralist}
\usepackage{enumitem} 
\setdefaultenum{(1)}{(a)}{(i)}{}
\setlist[enumerate,1]{label={\upshape(\arabic*)}}
\setlist[enumerate,2]{label={\upshape(\alph*)},ref=\theenumi\upshape(\alph*)}
\setlist[enumerate,3]{label={\upshape(\roman*)},ref=\theenumi\theenumii\upshape(\roman*)}
\usepackage[capitalize]{cleveref}
\crefname{Prop}{Proposition}{Propositions}
\crefname{Thm}{Theorem}{Theorems}
\crefname{Lem}{Lemma}{Lemmas}
\crefname{enumi}{Case}{Cases}

\newcommand{\prs}[1]{\left( #1 \right)}

\def\N{\ensuremath{\mathbb{N}}}
\def\P{\ensuremath{\mathbb{P}}}

\def\R{\ensuremath{\mathbb{R}}}

\def\Z{\ensuremath{\mathbb{Z}}}

\def\alg{\mathrm{alg}}
\def\Amp{\mathrm{Amp}}

\def\ch{\mathop{\mathrm{ch}}\nolimits}

\def\Coh{\mathop{\mathrm{Coh}}\nolimits}

\def\dim{\mathop{\mathrm{dim}}\nolimits}

\def\Ext{\mathop{\mathrm{Ext}}\nolimits}
\def\ext{\mathop{\mathrm{ext}}\nolimits}

\def\GL{\mathop{\mathrm{GL}}\nolimits}
\def\Hal{H^*_{\alg}}

\def\Hom{\mathop{\mathrm{Hom}}\nolimits}

\def\Nef{\mathrm{Nef}}

\def\NS{\mathop{\mathrm{NS}}\nolimits}

\def\Pic{\mathop{\mathrm{Pic}}\nolimits}

\def\rk{\mathop{\mathrm{rk}}}

\def\MG13{\ensuremath{{\mathcal M}_{\Gamma_1(3)}}}
\def\tildeMG13{\ensuremath{\widetilde{\mathcal M}_{\Gamma_1(3)}}}
\def\Stab{\mathop{\mathrm{Stab}}\nolimits}
\def\Stabd{\mathop{\Stab^{\dagger}}\nolimits}

\def\Db{\mathrm{D}^{\mathrm{b}}}

\newcommand\TFILTB[3]{%
\xymatrix@=1pc{
{0 = {#1}_0} \ar[rr]&&
{{#1}_1} \ar[rr]\ar[ld] &&
{{#1}_2} \ar[r]\ar[ld] &
{\cdots} \ar[r] & { {#1}_{#3-1}} \ar[rr] &&
{{#1}_{#3} = {#1}} \ar[ld]
\\
& *{{#2}_1} \ar@{.>}[ul] &&
{{#2}_2} \ar@{.>}[ul] & &&&
{{#2}_{{#3}}} \ar@{.>}[ul]
}}

\makeatletter
\newtheorem*{rep@theorem}{\rep@title}
\newcommand{\newreptheorem}[2]{%
\newenvironment{rep#1}[1]{%
 \def\rep@title{#2 \ref{##1}}%
 \begin{rep@theorem}}%
 {\end{rep@theorem}}}
\makeatother

\newtheorem{Thm}{Theorem}[section]
\newreptheorem{Thm}{Theorem}
\newtheorem{Prop}[Thm]{Proposition}

\newtheorem{Lem}[Thm]{Lemma}

\newtheorem{Cor}[Thm]{Corollary}
\newreptheorem{Cor}{Corollary}

\newreptheorem{Con}{Conjecture}

\newtheorem*{theorem*}{Theorem}
\newtheorem*{lemma*}{Lemma}
\newtheorem*{proposition*}{Proposition}
\newtheorem*{conjecture*}{Conjecture}
\newtheorem*{corollary*}{Corollary}
\newtheorem*{problem*}{Problem}

\newtheorem{Thm-int}{Theorem}

\theoremstyle{definition}
\newtheorem{Def-s}[Thm]{Definition}
\newtheorem{Def}[Thm]{Definition}
\newtheorem{Rem}[Thm]{Remark}

\newtheorem{Ex}[Thm]{Example}

\newenvironment{NB}{
\color{red}{\bf NB}. \footnotesize 
}{}

\def\N{\ensuremath{\mathbb{N}}}
\def\P{\ensuremath{\mathbb{P}}}

\def\R{\ensuremath{\mathbb{R}}}
\def\Z{\ensuremath{\mathbb{Z}}}

\def\BB{\ensuremath{\mathcal B}}
\def\CC{\ensuremath{\mathcal C}}

\def\GG{\ensuremath{\mathcal G}}

\def\MM{\ensuremath{\mathcal M}}

\def\OO{\ensuremath{\mathcal O}}
\def\PP{\ensuremath{\mathcal P}}

\def\QQ{\ensuremath{\mathcal Q}}

\def\WW{\ensuremath{\mathcal W}}

\def\u{\mathbf{u}}
\def\v{\mathbf{v}}
\def\w{\mathbf{w}}

\def\S{\mathbf{S}}

\newcommand{\leqor}{\underset{{\scriptscriptstyle (}-{\scriptscriptstyle )}}{<}}

\def\cal{\mathcal}
\def\Bbb{\mathbb}

\newcommand{\ignore}[1]{}
\begin{document}

\title{Weak Brill-Noether on abelian surfaces}
\author{Izzet Coskun}
\address{Department of Mathematics, Statistics, and CS \\
University of Illinois at Chicago, Chicago IL 60607}
\email{icoskun@uic.edu}
\author{Howard Nuer}
\address{Faculty of Mathematics \\
Technion, Israel Institute of Technology}
\email{hnuer@technion.ac.il}
\author{K\={o}ta Yoshioka}
\address{Department of Mathematics, Faculty of Science, Kobe University, Kobe, 657, Japan}
\email{yoshioka@math.kobe-u.ac.jp}

\thanks{During the preparation of this article, I.C. was partially supported by NSF grant DMS-2200684,  and K.Y. was partially supported by Grant in Aid for Scientific Research No. 21H04429, 23K03052 JSPS}

\begin{abstract}
We study the cohomology of a general stable sheaf on an abelian surface.
We say that a moduli space satisfies weak Brill-Noether if the general sheaf has at most one non-zero cohomology group.
Let $(X,H)$ be a polarized abelian surface and let $\v=(r, \xi, a)$ be a Mukai vector on $X$ with $\v^2 \geq 0$, $r>0$ and $\xi \cdot H>0$. We show that if $\rho(X)=1$ or $\rho(X)=2$ and $X$ contains an elliptic curve, then all the moduli spaces $M_{X,H}(\v)$ satisfy weak Brill-Noether. Conversely, if $\rho(X)>2$ or $\rho(X)=2$ and $X$ does not contain an elliptic curve, we show that there are infinitely many moduli spaces  $M_{X,H}(\v)$ that fail weak Brill-Noether. As a consequence, we classify Chern classes of Ulrich bundles on abelian surfaces.  

    \end{abstract}

\maketitle

\section{Introduction}
The Brill-Noether theory of projective curves has been a central pillar of algebraic geometry since the 19th century (see \cite{acgh}). On a smooth projective curve the cohomology of a general vector bundle of rank $r$ and degree $d$ is determined by its Euler characteristic. In contrast, on surfaces a general stable bundle with a fixed Chern character may have more than one nonzero cohomology group (see \cite{CoskunHuizengaNuer, CoskunNuerYoshioka} for examples). Determining the cohomology of a general stable bundle on a higher dimensional variety can be a challenging problem. In this paper, we address the problem for abelian surfaces.

\smallskip

Let $(X,H)$ be a polarized abelian surface and let $M_{X,H}({\bf v})$ denote the moduli space of $H$-Gieseker semistable sheaves on $X$ with Mukai vector ${\bf v} = (r, \xi, a)$. Given two Mukai vectors $\v, \v'$, let $\langle \v, \v' \rangle$ denote their Mukai pairing. Denote the Picard rank of $X$ by  $\rho(X)$, which is an integer between 1 and 4.  Let $d(\NS(X))$ denote the discriminant of the N\'{e}ron-Severi lattice.

\begin{Def}
    A moduli space  $M_{X, H}({\bf v})$ satisfies the {\em weak Brill-Noether property} if there exists a sheaf $E \in M_{X, H}({\bf v})$ such that $E$ has at most one nonzero cohomology group. 
\end{Def} 

\begin{Rem}\label{Rem:WBNDetermingsCohomology}
    When $M_{X, H}(\v)$ is irreducible and satisfies weak Brill-Noether, then the general sheaf $E$ in $M_{X, H}$ has at most one nonzero cohomology group. Furthermore, the cohomology of $E$ is determined by numerical data:
    \begin{enumerate}
        \item  If the Euler characteristic $\chi({\bf v}) \leq 0$, then  $h^1(E) = - \chi({\bf v})$ and $h^0(E)=h^2(E)=0$.
        \item If $\chi({\bf v}) > 0$ and $\xi \cdot H >0$, then $h^0(E)= \chi({\bf v})$ and $h^1(E)=h^2(E)=0$.
         \item If instead $\chi({\bf v}) > 0$ and $\xi \cdot H <0$, then $h^2(E)= \chi({\bf v})$ and $h^0(E)=h^1(E)=0$.
          
    \end{enumerate}
 In particular, if $M_{X,H}({\bf v})$ satisfies weak Brill-Noether and $\chi({\bf v}) = 0$, then $h^i(E)=0$ for all $i$. Hence, the weak Brill-Noether property plays a central role in the construction of theta divisors and Ulrich bundles.    
\end{Rem}

Abelian surfaces form a natural class of surfaces on which to study the weak Brill-Noether problem.  On $\P^2$, G\"{o}ttsche and Hirschowitz proved that the general stable bundle  satisfies weak Brill-Noether \cite{GottscheHirschowitz}. When a surface $X$ contains curves with negative self-intersection, then there are moduli spaces of stable sheaves on $X$ for which weak Brill-Noether fails \cite[Proposition 5.10]{CoskunHuizengaNuer}. On surfaces such as Hirzebruch surfaces and certain del Pezzo surfaces, the negative self-intersection curves fully account for the failure of weak Brill-Noether (see \cite{CoskunHuizengaWBN, CoskunHuizengaBN, LevineZhang}). On more general blow-ups of $\P^2$ the situation is more complicated (see \cite{CoskunHuizenga23}).  Like $\P^2$, Abelian surfaces do not contain curves of negative self-intersection. 

The higher cohomology of line bundles is another source of failure of weak Brill-Noether. For example, on K3 surfaces of Picard rank 1, the higher cohomology of the structure sheaf yields many examples of moduli spaces that fail weak Brill-Noether (see \cite{CoskunNuerYoshioka}). For abelian surfaces of Picard rank 1, the general line bundle has expected cohomology. Hence, one can expect weak Brill-Noether to hold on abelian surfaces of Picard rank 1. Indeed, our main theorem in this paper is the following. 

\begin{Thm}\label{Thm:MainTheorem1}
Let $(X,H)$ be a polarized abelian surface. Let ${\bf v}= (r, \xi, a)$ be a Mukai vector with $\v^2 \geq 0$, $r >0$, and $\xi \cdot H >0$. 
\begin{enumerate}
    \item If $\rho(X)=1$, then weak Brill-Noether holds for  $M_{X, H}(\v)$.
    \item If $\rho(X) =2$ and $-d(\NS(X))$ is a square, then weak Brill-Noether holds for $M_{X, H}(\v)$. 
\end{enumerate}
Conversely, if $\rho(X) \geq 3$ or if $\rho(X)=2$ and $-d(\NS(X))$ is not a square, then there exists infinitely many $\v$ with $\v^2\geq0$ and $r>0$ such that for some $H$ with $\xi \cdot H>0$, the moduli space $M_{X, H}(\v)$ does not satisfy weak Brill-Noether.
\end{Thm}

\begin{Rem}
  The third author previously obtained  the case $\rho(X)=1$ in \cite[Lemma 3.13]{Y:Stability} by different methods.   
\end{Rem}
  
\begin{Rem}
    The weak Brill-Noether property for $M_{X, H}(\v)$ for $\v$ with $\v^2 \geq 0$ and $r>0$ such that $\xi \cdot H < 0$ can be studied using Serre duality (see Section \ref{sec:NegativeSlope}).
\end{Rem}

We will prove Theorem \ref{Thm:MainTheorem1} using Fourier-Mukai transforms and Bridgeland stability and provide detailed information about Mukai vectors $\v$ for which $M_{X,H}(\v)$ does not satisfy weak Brill-Noether.  To state our main technical theorem we need the following definition.

\begin{Def}\label{Def:IsotropicDecomposition}
    A Mukai vector $\v$ admits an {\em isotopric decomposition} if $\v= \ell_1 \v_1 + \ell_2 \v_2$, where $\v_i = (r_i, \xi_i, a_i)$ with $$r_i >0, \quad (\xi_i \cdot H)(\xi\cdot H) >0, \quad a_1 > 0, \quad a_2<0,$$ 
$$\v_1^2=\v_2^2 = 0, \quad \langle\v_1, \v_2 \rangle=1, \quad \mbox{and} \quad \{ \ell_1,\ell_2 \}=\left\{ \frac{\v^2}{2},1\right\}.$$
\end{Def}

\begin{Thm}[Theorem \ref{thm:wbn}] \label{Thm:introtech}
    Let $(X,H)$ be a polarized abelian surface. Let  $\v=(r,\xi,a)$ be a Mukai vector such $\v^2 \geq 0$,  
$r>0$, and $\xi \cdot H >0$. If weak Brill-Noether does not hold for $M_{X,H}(\v)$, then $\v$ admits an isotropic decomposition.
\end{Thm}

  We obtain the following useful corollaries of independent interest.

\begin{Cor}[Corollary \ref{Cor:weak}]\label{Cor:introweak}
   Let $(X,H)$ be a polarized abelian surface. Let $\v= (r, \xi, a)$ be a Mukai vector such that $\v^2 \geq 0$, $r>0$, and $\xi \cdot H > 0$.
   \begin{enumerate}
       \item If $\v$ is not primitive, then weak Brill-Noether holds for $M_{X, H}(\v)$.
       \item If $\v$ is primitive and $\v^2 + 2 > 2r$, then weak Brill-Noether holds for $M_{X, H}(\v)$.
       \item If $\v_j = (r, \xi, a-j)$ for an integer $j >0$, then weak Brill-Noether holds for $M_{X, H}(\v_j)$.
   \end{enumerate}
\end{Cor}

\begin{Rem}
It is important to mention that  \cref{Thm:introtech,Cor:introweak} completely solve the weak Brill-Noether problem for abelian surfaces.  Weak Brill-Noether always holds for $M_{X,H}(\v)$ except possibly for a primitive Mukai vector $\v=(r,\xi,a)$ with $a=\lfloor\frac{\xi^2}{2r}\rfloor$.  In this case, we simply check by hand if $\v$ admits an isotropic decomposition, in which case we compute the cohomology of the general sheaf in \cref{thm:weak}. Otherwise, weak Brill-Noether holds in this case as well. 
\end{Rem}

As a consequence of Theorem \ref{thm:wbn}, we obtain a classification of Ulrich bundles on abelian surfaces. Previously, Beauville \cite{Be} constructed Ulrich bundles of rank 2 on abelian surfaces.

\begin{Thm}[Theorem \ref{thm-ulrich}]
    Let $(X,H)$ be a polarized abelian surface. If $E$ is an Ulrich bundle of rank $r\geq 2$, then its Mukai vector $\v$ has the form $(r, \xi,0)e^H$ with $2\xi \cdot H = rH^2$ and $\xi^2 \geq 0$. Conversely, suppose that $\v$ is a Mukai vector satisfying these properties and $(r, \xi, 0)$ does not admit an isotropic decomposition, then the general bundle in $M_{X,H}(\v)$ is an Ulrich bundle.
\end{Thm}

\subsection*{Organization of the paper} In Section \ref{sec-prelim}, we introduce the necessary background on abelian surfaces, moduli spaces of sheaves on abelian surfaces, Fourier-Mukai transforms and Bridgeland stability. In Section \ref{sec-cohomology}, we prove Theorem \ref{thm:wbn} and obtain criteria for weak Brill-Noether to hold. In Section \ref{sec-Neron-Severi}, we determine when a Mukai vector can have an isotropic decomposition and prove Theorem \ref{Thm:MainTheorem1}. In Section \ref{sec-Ulrich}, we give applications to classification of Ulrich bundles on abelian surfaces.

\subsection*{Acknowledgments} We would like to thank Arend Bayer, Aaron Bertram, Jack Huizenga, Jesse Kass, Emanuele Macr\`i, Ari Shnidman, and Paolo Stellari for many invaluable discussions.

\section{Background results}\label{sec-prelim}
In this section, we recall the necessary background results on moduli spaces of sheaves on abelian surfaces, Fourier-Mukai transforms and Bridgeland stability.

\subsection{Mukai Lattice}
Let $X$ be an abelian surface. Then the even cohomology $H^{2*}(X.{\Bbb Z})$ of $X$
$$H^{2*}(X,{\Bbb Z}):={\Bbb Z} \oplus H^2(X,{\Bbb Z}) \oplus {\Bbb Z}$$ has a natural pairing called  the {\em Mukai pairing} defined by 
\begin{equation*}
\langle x,y \rangle:=x_1 \cdot y_1-x_0 y_2-y_0 x_2 \in {\Bbb Z}
\end{equation*}
for $x=(x_0,x_1,x_2), y=(y_0,y_1,y_2) \in  {\Bbb Z} \oplus H^2(X,{\Bbb Z}) \oplus {\Bbb Z}$, where $x_1 \cdot y_1$ is the cup product on $H^2(X, {\Bbb Z})$.
The pair $\prs{H^{2*}(X,{\Bbb Z}),\langle\;\;,\;\; \rangle}$ is called the {\em Mukai lattice} of $X$. For simplicity, when we take the self-product of a divisor or a Mukai vector we drop the brackets.

The Mukai lattice has a natural Hodge structure. Denoting the N\'{e}ron-Severi space of $X$ by $\NS(X)$, the algebraic part of $H^{2*}(X,{\Bbb Z})$ is 
$$H^{2*}(X,{\Bbb Z})_{\alg}:={\Bbb Z} \oplus \NS(X) \oplus {\Bbb Z}.$$ 

For an object $E$ of the bounded derived category $\Db(X)$ of the category $\Coh(X)$ of
coherent sheaves, 
$\v(E):=\ch(E) \in H^{2*}(X,{\Bbb Z})$ is the Mukai vector of $E$. 
We also say that $\v \in H^{2*}(X,{\Bbb Z})_{\alg}$ is a Mukai vector.  
By the Hirzebruch-Riemann-Roch Theorem, the Mukai pairing is defined so that for any $E,F\in\Db(X)$, we have 
$$\langle\v(E),\v(F)\rangle=-\chi(E,F)=-\sum_i(-1)^i\dim\Ext^i(E,F),$$
where $\chi\prs{\;\;,\;\;}$ is the Euler pairing.

\subsection{Stable sheaves} Let $E$ be a coherent sheaf of pure dimension $d$ on $X$. Then the Hilbert and the reduced Hilbert polynomials $P_{E,H}(m)$ and $p_{E,H}(m)$, respectively, of $E$ with respect to an ample $H$ on $X$ are defined by 
$$P_{E,H}(m)= \chi(E(mH))= a_d \frac{m^d}{d!} + \mbox{l.o.t.},  \quad p_{E,H}(m) = \frac{P_{E,H}(m)}{a_d}.$$ The sheaf $E$ is called $H$-Gieseker (semi)stable if every proper nonzero subsheaf $F$ satisfies  $$p_{F,H}(m) \leqor p_{E,H}(m)$$
for $m \gg 0$. When the ample $H$ is clear from context, we simply say the sheaf is (semi)stable. When $E$ is torsion free, let  $$\mu_H(E) := \frac{c_1(E) \cdot H}{\rk(E)}.$$ The sheaf $E$ is $\mu_H$-(semi)stable if for every subsheaf $0 \ne F \subsetneq E$ of smaller rank, we have $$\mu_H(F) \leqor \mu_H(E).$$

If $E$ is a stable sheaf on an abelian surface $X$, then 
\cite[Prop. 3.8]{Mu:Fourier} implies $\v^2(E) \geq 0$.

Gieseker \cite{Gieseker} and Maruyama \cite{Maruyama} constructed projective moduli spaces parameterizing $S$-equivalence classes of semistable sheaves on any projective variety. In the case of abelian surfaces, these moduli spaces have been intensively studied  by Mukai, Yoshioka and others. 
Let $M_{X,H}(\v)$ denote the moduli space of $S$-equivalence classes of $H$-Gieseker semistable sheaves of Mukai vector $\v$. 
Denote by $\MM_{X,H}(\v)$ the stack of $H$-Gieseker semistable sheaves of Mukai vector $\v$ and by $\MM_{X,H}(\v)^{\mu ss}$ the stack of $\mu_H$-semistable sheaves of Mukai vector $\v$.  When $X$ and $H$ are understood from context, we may drop them from the notation. We now summarize basic facts about these moduli spaces.

\begin{Thm}\label{Thm:ClassicalFactsModuli}
    Let $(X,H)$ be a polarized abelian surface. Let $\v$ be a Mukai vector such that $\v^2 \geq 0$. Assume that $\v= m \v_0$ for a primitive Mukai vector $\v_0$. 
    \begin{enumerate}
        \item\cite[Theorem 3.18]{Y:TwistedFM1} If $\v^2 >0$, then $M_{X,H}(\v)$ is irreducible of dimension $\v^2 +2$. 
        \item\cite[Proposition 4.2]{Y:abel} If $\v^2 =0$, then $M_{X,H}(\v)$ is irreducible of dimension $2m$. The general sheaf in $M_{X,H}(\v)$ is $\mu_H$-stable if and only if $m=1$.  When $m>1$ every $E\in M_{X,H}(\v)$ is $S$-equivalent to $\bigoplus_{i=1}^m E_i$ with $E_i\in M_{X,H}(\v_0)$.  Moreover, in this case, stability and semi-stability do not depend on the choice of $H$.
    \end{enumerate}
Furthermore, if $\v$ is primitive, then $M_{X,H}(\v)$ is smooth.
\end{Thm}

\begin{Lem}[{\cite[Proposition 3.5]{KY}}]\label{lem:mu-s}
Let $(X,H)$ be a polarized abelian surface. Assume that $r \geq 2$, $\v^2>0$, and
$\v \ne (r,0,-1) e^{\eta}$ $(\eta \in \NS(X))$.
Then for a general $H$ there is a $\mu_H$-stable vector bundle with $\v(E)=\v$.
\end{Lem}
\begin{Lem}[{\cite[Lemma 3.8]{KY}}]\label{Lem:ReducibleModuli}
Let $(X,H)$ be a polarized abelian surface. 
Assume that $H$ belongs to a wall with respect to $\v=(r,\xi,a)$ and $H'$ is a general ample divisor which is close to $H$.  For the stack $\MM_H(\v)^{\mu ss}$ 
of $\mu_H$-semistable sheaves,
$$\dim \MM_H(\v)^{\mu ss}=\v^2+1.$$  If $\MM_H(\v)^{\mu ss}\backslash \MM_{H'}(\v)^{\mu ss}$ contains a $(\v^2+1)$-dimensional component, then there is a Mukai vector $\w=(r',\xi',a')$ such that $r'>0$, $(r\xi'-r'\xi)\cdot H=0$, $\w^2=0$ and $\langle \v,\w\rangle=1$.  A general element of $\MM_H(\v)^{\mu ss}\backslash \MM_{H'}(\v)^{\mu ss}$ fits in an exact sequence 
$$0\to E_1\to E\to E_2\to 0$$
where $E_i\in \MM_{H'}(\v_i)^{ss}$, $\{\v_1,\v_2\}=\{l'\w,l''\w+\u\}$, $l',l''\ge 0$, $\u^2=0$ and $\langle \u,\w\rangle=1$.
\end{Lem}

\subsection{Semi-homogeneous bundles} We begin by studying the cohomology of vector bundles in the extremal case $\v^2(E)=0$. These are the semi-homogeneous bundles on $X$. These bundles will play a central role in our general analysis, as we will explain in the next subsection.
\begin{Def}
    Let $t_x$ denote translation by $x \in X$. A nonzero torsion free sheaf $E$ on $X$ is called {\em semi-homogeneous} if for every point $x \in X$, there exists a line bundle $L_x$ such that $t_x^* E \cong E \otimes L_x$.  In particular, $E$ must be a vector bundle. 
\end{Def}
The following is a useful characterization of semi-homogeneous bundles that is due to Mukai, which we state in the case of dimension two:
\begin{Prop}[{\cite[Theorem 5.8]{Mu}}]\label{Prop:SemiHomogEquivalence}
    Let $E$ be a simple torsion free sheaf on an abelian surface $X$.  Then the following conditions are equivalent:
    \begin{enumerate}
        \item $\ext^1(E,E)=2$;
        \item $E$ is semi-homogeneous;
        \item there exists an isogeny $\pi\colon Y\to X$ and a line bundle $L$ on $Y$ such that $E\cong\pi_*(L)$.
    \end{enumerate}
    In particular, for a simple torsion free sheaf $E$, the sheaf $E$ is semi-homogeneous if and only if $\v(E)^2=0$.
\end{Prop}

Let $\widehat{X}$ be the dual abelian surface of $X$ and
${\cal P}$ be the Poincar\'{e} line bundle on $X \times \widehat{X}$. Then following proposition characterizes the cohomology of semi-homogeneous bundles.

\begin{Prop}\label{Prop:CohomologySemiHomogeneous}
    Let $E$ be a simple semi-homogeneous vector bundle  on an abelian surface $X$ with $\v(E)=(r,\xi,a)$.  
    \begin{enumerate}
\item[(1)]
Assume that $(\xi^2)>0$. 
\begin{enumerate}
\item
If $\xi \cdot H>0$, then for all $t\in\widehat{X}$, $E\otimes\PP_t$ has $H^0$ while all other cohomology vanishes.
\item
If $\xi \cdot H<0$, then for all $t\in\widehat{X}$, $E\otimes\PP_t$ has $H^2$ while all other cohomology vanishes.
\end{enumerate}
\item[(2)]
Assume that $(\xi^2)<0$, then for all $t\in\widehat{X}$, $E\otimes\PP_t$ has $H^1$ while all other cohomology vanishes.
\item[(3)]
Assume that $(\xi^2)=0$ with $\xi \ne 0$.  Then $a=0$.
\begin{enumerate}
    \item If $\xi \cdot H>0$, then for general $t\in\widehat{X}$, $E\otimes\PP_t$ has no cohomology, while there exists some $t\in\widehat{X}$ such that $E\otimes\PP_t$ has both $H^0$ and $H^1$.
    \item If $\xi \cdot H<0$, then for general $t\in\widehat{X}$, $E\otimes\PP_t$ has no cohomology, while there exists some $t\in\widehat{X}$ such that $E\otimes\PP_t$ has both $H^1$ and $H^2$.
\end{enumerate}
\end{enumerate}
\end{Prop}
\begin{proof}
By \cref{Prop:SemiHomogEquivalence}, there is an isogeny $\pi:Y \to X$ and a line bundle $L$ on $Y$
such that $\pi_*(L)=E$.
If $(\xi^2)>0$, then $(c_1(L)^2) >0$ so that
$L$ is ample or $L^{\vee}$ is ample.
If $L$ is ample, then $0 \ne H^0(Y,L)=H^0(X,E)$.
If $L^{\vee}$ is ample, then $0 \ne H^2(Y,L)=H^2(X,E)$.
In both cases, all other cohomology vanishes.
If $(\xi^2)<0$, then $(c_1(L)^2)<0$, so 
$0 \ne H^1(Y,L)=H^1(X,E)$, while all other cohomology vanishes.
Assume that $(\xi^2)=(c_1(L)^2)=0$.
Then $c_1(L)$ is nef or $-c_1(L)$ is nef.
If $c_1(L)$ is nef, then
there is $t \in \widehat{X}$ such that
$0 \ne H^1(Y,L \otimes \pi^*({\cal P}_t))=H^1(X,E \otimes {\cal P}_t)$ and the same is true for $H^0(X,E\otimes\PP_t)$.
If $-c_1(L)$ is nef, then similarly there is $t \in \widehat{X}$ such that
$0 \ne H^2(Y,L \otimes \pi^*({\cal P}_t))=H^2(X,F \otimes {\cal P}_t)$ and the same is true for $H^1(X,E\otimes\PP_t)$.
\end{proof}

\begin{Rem}\label{Rem:CohomologySemiHomogeneous}
  Let $E$ be a semistable torsion free sheaf with $\v(E)^2=0$.  
  By \cref{Thm:ClassicalFactsModuli}, all of the factors in the Jordan-H\"{o}lder filtration of $E$ are simple semi-homogeneous bundles.  Hence, \cref{Prop:CohomologySemiHomogeneous} computes the cohomology of those bundles as well.   
\end{Rem}

\subsection{Fourier-Mukai transforms} 

Let $\Phi_{X \to \widehat{X}}^{{\cal P}^{\vee}}:\Db(X) \to \Db(\widehat{X})$
be the Fourier-Mukai transform whose kernel is ${\cal P}^{\vee}$ \cite{Mukai:1981}:
\begin{equation}
\begin{matrix}
\Phi_{X \to \widehat{X}}^{{\cal P}^{\vee}}: & \Db(X) & \to & \Db(\widehat{X})\\
& E & \mapsto & {\bf R}p_{\widehat{X}*}({\cal P}^{\vee} \otimes p_X^*(E)),
\end{matrix}
\end{equation}
where $p_X:X \times \widehat{X} \to X$ and
$p_{\widehat{X}}:X \times \widehat{X} \to \widehat{X}$
are the two projections.
We set
$\Phi:=\Phi_{X \to \widehat{X}}^{{\cal P}^{\vee}}$ and
$\widehat{\Phi}:=\Phi_{\widehat{X} \to X}^{{\cal P}}$.
We also set 
\begin{equation}
\begin{split}
\Phi^i(E):=& H^i(\Phi(E)) \in \Coh(\widehat{X}),\; E \in \Coh(X)\\
\widehat{\Phi}^i(F):=& H^i(\widehat{\Phi}(F)) \in \Coh(X),\;F \in \Coh(\widehat{X}).
\end{split}
\end{equation}

The Fourier-Mukai transform $\Phi$ induces an isometry of Hodge structures
\begin{equation}
\begin{matrix}
H^{2*}(X,{\Bbb Z}) & \to & H^{2*}(\widehat{X},{\Bbb Z})\\
(r,\xi,a) & \mapsto & (a,-\widehat{\xi},r),
\end{matrix}
\end{equation}
where $\widehat{\xi}$ is the Poincar\'{e} dual of $\xi$.
If $\Phi^i(E)$ is a torsion sheaf, then by cohomology and base change it follows that
$$H^i(X,E \otimes {\cal P}^{\vee}_{|X \times \{ y \}}) \otimes k(y) \cong
\Phi^i(E) \otimes k(y)=0$$
for a general $y \in \widehat{X}$.
Hence the Fourier-Mukai transform can be used to study
cohomology groups of a general stable sheaf.
More generally, the connection between the Fourier-Mukai transform and studying the cohomology groups of a general stable sheaf is provided by the following result due to Mukai.
\begin{Lem}\label{Lem:FundamentalLemmaAbelian}
Let $X$ be an abelian surface and $E\in\Coh(X)$.  If $\hat{E}=\Phi(E)[i]\in\Coh(\widehat{X})$, then 
$$H^j(X,E\otimes\PP_t)=\Ext_{\hat{X}}^{j+2-i}(\OO_{-t},\hat{E})$$
and
$$\Ext_X^j(\OO_x,E)\cong H^{j-i}(\hat{X},\hat{E}\otimes\PP_{x})$$
for every $x\in X$, $t\in\hat{X}$, and $j\geq 0$, .
\end{Lem}
\begin{proof}
This is proved as in \cite[Prop. 2.7]{Mukai:1981} except that Mukai uses the transform with $\Phi$ with kernel $\PP$ instead of $\PP^\vee$.
\end{proof}

In particular, if $\Phi(E)\in\Coh(\widehat{X})$, then $H^j(X,E\otimes\PP_t)=\Ext^{j+2}(\OO_{-t},\hat{E})=0$ for any $j>0$ and any $t\in\hat{X}$, and thus $E$ has at most one nonzero cohomology group.

The following follows from the structure of Fourier-Mukai transforms
on an abelian surface \cite{BridgelandMacioca,Or1} or \cite[Prop. 4.4]{Y:abel}.
\begin{Lem}\label{lem:semihom}
Let $E$ be a stable sheaf with $\v(E)=(r,\xi,a)$ and $\v^2=0$.
\begin{enumerate}
\item[(1)]
Assume that $(\xi^2)>0$. 
\begin{enumerate}
\item
If $\xi \cdot H>0$, then $\Phi(E) \in \Coh(\widehat{X})$. 
\item
If $\xi \cdot H<0$, then $\Phi(E)[2] \in \Coh(\widehat{X})$.
\end{enumerate}
\item[(2)]
Assume that $(\xi^2)<0$. Then $\Phi(E)[1] \in \Coh(\widehat{X})$. 
\item[(3)]
Assume that $(\xi^2)=0$ with $\xi \ne 0$.
\begin{enumerate}
\item
If $r=0$, then $\xi$ is effective.
Hence $\Phi(E) \in \Coh(\widehat{X})$ for $a> 0$
and $\Phi(E)[1] \in \Coh(\widehat{X})$ for $a \leq 0$.
\item 
If $r>0$, then $a=0$. If $\xi \cdot H>0$, then
$\Phi(E)[1] \in \Coh(\widehat{X})$. 
If $\xi \cdot H<0$, then
$\Phi(E)[2] \in \Coh(\widehat{X})$. 
\end{enumerate} 
\end{enumerate}
\end{Lem}

\begin{proof}
Assume that $(\xi^2) \ne 0$.
As $r>0$, the claims follow from Proposition \ref{Prop:CohomologySemiHomogeneous} and 
Remark \ref{Rem:CohomologySemiHomogeneous} using cohomology and base change.

We next assume that $(\xi^2)=0$ and $\xi \ne 0$.
Then $ra=0$ and $\xi$ is nef or $-\xi$ is nef, according to the sign of $\xi\cdot H$.
Assume first that $r=0$. Then $\xi$ is effective and thus $\xi\cdot H>0$.
Hence as above it follows from the stability of $E$ that $\Phi^2(E)=0$.  
Thus $\Phi(E) \in \Coh(\widehat{X})$ for $a>0$ and
$\Phi(E)[1] \in \Coh(\widehat{X})$ for $a \leq 0$.

If $r>0$, then  $ra=0$ implies $a=0$.
If $\xi \cdot H>0$, then the torsion freeness of $\Phi^0(E)$
and Proposition \ref{Prop:CohomologySemiHomogeneous} (3) (a) imply
$\Phi^0(E)=\Phi^2(E)=0$.
Thus $\Phi(E)[1] \in \Coh(\widehat{X})$.
If $\xi \cdot H<0$, then 
$H^0(X,E \otimes {\cal P}_t^{\vee})=0$ for all $t \in \widehat{X}$.
By the base change theorem, $\Phi^0(E)=0$ and
$\Phi^1(E)$ is torsion free. 
By Propoisition \ref{Prop:CohomologySemiHomogeneous} (3) (b),
$\Phi^1(E)=0$. Thus $\Phi(E)[2] \in \Coh(\widehat{X})$.

\end{proof}

Having resolved the question of the weak Brill-Noether property for stable sheaves with isotropic Mukai vector, and thus all semi-homogeneous bundles, we will go on to address the weak Brill-Noether property for semistable sheaves of positive-square Mukai vector.  We quote the following result which is used to prove our main result.
\begin{Thm}[{\cite[Thm. 3.1]{Y:abel}}]\label{thm:abel}
Assume that $\v$ is primitive.
For a general $E \in {\cal M}_H(\v)$,
assume that $\Phi(E) \not \in \Coh(\widehat{X})$.
If $\Phi^i(E)=0$ for $i>1$, then $\Phi^0(E)$ and $\Phi^1(E)$ are semistable sheaves with
$\v(\Phi^i(E))=\ell_i \w_i$, $\w_i^2 =0$, $\langle \w_0,\w_1 \rangle=-1$ and $(\ell_0,\ell_1)=(\ell,1),(1,\ell)$,
where $\ell=\langle \v^2 \rangle/2$.
\end{Thm}

\subsection{Stability conditions on an abelian surface}
We briefly recall Bridgeland  stability conditions on an abelian surface $X$.
For more details, see \cite{Bridgeland}, \cite{BridgelandK3} and \cite{BridgelandSSC}.
A stability condition $\sigma=(Z_\sigma,\QQ_\sigma)$ 
on $\Db(X)$
consists of a group homomorphism $Z_\sigma: \Db(X) \to {\Bbb C}$
and a slicing $\QQ_\sigma$ of $\Db(X)$ such that if
$0 \ne E \in \QQ_\sigma(\phi)$, then 
$Z_\sigma(E) =m(E)\exp(\pi \sqrt{-1} \phi)$ 
for some $m(E) \in {\Bbb R}_{>0}$.
We further require  that $\sigma$ satisfies certain axioms such as the support property, but since we do not use these explicitly we refer the reader to the references above for the definitions.

The set of stability conditions has the structure of a complex 
manifold.
We denote this space by $\Stab(X)$.

\begin{Def}\label{defn:B-stable}
\begin{enumerate}
\item[(1)]
A $\sigma$-semistable object $E$ of phase $\phi$ is
an object of  $\QQ_\sigma(\phi)$.
If $E$ is a simple object of $\QQ_\sigma(\phi)$, then
$E$ is $\sigma$-stable. 
\item[(2)]
For $\v \in H^*(X,{\Bbb Z})_{\alg}$,
we denote the moduli stack of $\sigma$-semi-stable objects $E$ with 
$\v(E)=\v$ by ${\cal M}_\sigma(\v)$,
where we usually choose 
$\phi:=\frac{\mathrm{Im}(\log Z_\sigma(v))}{\pi} \in (-1,1]$.
\end{enumerate}
\end{Def}

Let $\Stab(X)$ be the space of stability conditions.
For an equivalence $\Phi:\Db(X) \to \Db(X')$,
we have an isomorphism $\Phi:\Stab(X) \to \Stab(X')$
such that 
$\Phi(\sigma)$ $(\sigma \in \Stab(X))$ 
is a stability condition given by
\begin{equation}\label{eq:FM-action}
\begin{split}
Z_{\Phi(\sigma)}=& Z_\sigma \circ \Phi^{-1}:\Db(X') \to {\Bbb C},\\
{\cal Q}_{\Phi(\sigma)}(\phi)=& \Phi(\QQ_\sigma(\phi)).
\end{split}
\end{equation}
We also have an action of the universal covering 
$\widetilde{\GL}_2^+({\Bbb R})$ of $\GL_2^+({\Bbb R})$ on 
$\Stab(X)$.
Since ${\Bbb C}^{\times} \subset  \GL_2^+({\Bbb R})$,
we have an injective homomorphism
${\Bbb C} \to  \widetilde{\GL}_2^+({\Bbb R})$.
Thus 
we have an action of $\lambda \in {\Bbb C}$ on
$\Stab(X)$. 
For a stability condition $\sigma \in \Stab(X)$, 
$\lambda(\sigma)$ is 
given by
\begin{equation}\label{eq:C-action}
\begin{split}
Z_{\lambda(\sigma)}=&
\exp(-\pi  \sqrt{-1} \lambda)Z_\sigma\\
{\cal Q}_{\lambda(\sigma)}(\phi)=&
\QQ_\sigma(\phi+\mathrm{Re}\lambda).
\end{split}
\end{equation}

\begin{Def}\label{defm:moduli}
Let ${\cal M}_{\sigma}(\v)$ be the moduli stack of
$\sigma$-stable objects $E$ with $\v(E)=\v$. 
\end{Def}


\subsection{Wall-and-chamber decomposition}
For a fixed Mukai vector $\v \in H^{2*}(X,{\Bbb Z})$, there exists a locally finite set of \emph{walls} (real codimension one submanifolds with boundary) in ${\cal H}$, 
depending only on $\v$, with the following properties (see \cite{BridgelandSSC}):
\begin{enumerate}
\item[(1)]
When $\sigma$ varies in a chamber, that is, a connected component of 
the complement of the union of walls, the sets of 
$\sigma$-semistable and $\sigma$-stable objects of class $\v$ do not change. 
If $\v$ is primitive, then $\sigma$-stability coincides with $\sigma$-semistability 
for $\sigma$ in a chamber for $\v$.

\item[(2)]
 When $\sigma$ lies on a wall $W \subset {\cal H}$, 
there is a $\sigma$-semistable object of class $\v$ that is unstable in one of the adjacent chambers 
and semistable in the other adjacent chamber.
 If $\sigma= (Z_\sigma, {\cal Q}_\sigma)$ lies on a wall, there exists a $\sigma$-semistable object $E$ 
of Mukai vector $\v$ and phase $\phi$ and a subobject $F \subset E$ in ${\cal Q}_\sigma(\phi)$ but 
$v(F) \not \in {\Bbb R} v$.

\item[(3)]
 Given a polarization $H\in\Amp(X)$ and the Mukai vector $\v$ of an $H$-Gieseker semistable sheaf, 
there exists a chamber ${\cal C}$ for $\v$, the \emph{Gieseker chamber}, 
where the set of $\sigma$-semistable objects of class 
$\v$ coincides with the set of $H$-Gieseker semistable sheaves \cite[Prop. 14.2]{BridgelandK3}.
\end{enumerate}
We set
$$
{\cal H}:=\NS(X)_{\Bbb R} \times \Amp(X)_{\Bbb R}.
$$
For $(\beta,\omega) \in {\cal H}$,
Bridgeland \cite{BridgelandSSC} constructed a stability condition
$\sigma_{(\beta,\omega)}=(Z_{(\beta,\omega)},{\cal A}_{(\beta,\omega)})$
which is characterized
by the stability of $k_x$ ($x \in X$).

\begin{Def}\label{defn:W}
For a Mukai vector $\v_1$,
let $W_{\v_1}$ be a closed subset of ${\cal H}$ such that
$$
W_{\v_1}=\{(\beta,\omega) \in {\cal H} \mid {\Bbb R} Z_{\sigma_{(\beta,\omega)}}(\v)=
 {\Bbb R}Z_{\sigma_{(\beta,\omega)}}(\v_1) \}.
$$
\end{Def}

For a wall $W$ in ${\cal H}$,
there is $\v_1$ such that $W= W_{\v_1}$.  
We note that a Mukai vector $\v_1$ defines a wall in the space of stability conditions
if and only if
\begin{equation}\label{eq:def-wall}
\langle \v_1^2 \rangle \geq 0,\quad \langle (\v-\v_1)^2 \rangle \geq 0,\quad \langle \v_1,\v-\v_1 \rangle>0,\quad 
\langle \v_1^2 \rangle \  \langle \v^2 \rangle < \langle \v_1,\v \rangle^2
\end{equation}
(see \cite[Def. 1.2 and Prop. 1.3]{Y}).
By \cite[Thm. 3.8, Prop. 5.14]{YY2}, we have the following result.
\begin{Prop}\label{prop:isom}
\begin{enumerate}
\item[(1)]
$\Phi[1]$ induces an isomorphism
$$
{\cal M}_{(0,tH)}(r,\xi,a) \to {\cal M}_{(0,t^{-1} \widehat{H})}(-a,\widehat{\xi},-r).
$$
\item[(2)]
Assume that $t \ll 1$.
\begin{enumerate}
\item
If $a \leq 0$, then
${\cal M}_{(0,t^{-1} \widehat{H})}(-a,\widehat{\xi},-r)={\cal M}_{\widehat{H}}(-a,\widehat{\xi},-r)$.
\item
If $a>0$, then
${\cal M}_{(0,t^{-1} \widehat{H})}(-a,\widehat{\xi},-r)$ 
consists of $E^{\vee}[1]$, $E \in {\cal M}_{\widehat{H}}(a,-\widehat{\xi},r)$.
\end{enumerate}
\end{enumerate}
\end{Prop}

\subsubsection{Totally semistable walls}

\begin{Def}
A wall $\WW$ is totally semistable if
$\MM_{\sigma_+}(\v) \cap \MM_{\sigma_-}(\v) =\emptyset$,
where $\sigma_\pm$ are in the adjacent two chambers. 

\end{Def}

We set
$$
{\cal I}_1:=\{ \v_1 \mid \v_1^2 =0, \langle \v,\v_1 \rangle=1 \}.
$$

\begin{Prop}[{\cite[Thm. 5.3.5]{MYY:2018}, see also \cite[Lem. 5.3.4]{MYY:2018} and \cite[Thm. 4.2]{NuerYoshioka}}]\label{prop:tot-wall}
Let $\WW$ be a wall for $\v$ with $\v^2\geq 0$.  Then $\WW$ is totally semistable if and only if $\WW=\WW_{\v_1}$ for some $\v_1 \in {\cal I}_1$.
In particular there is no totally semistable walls if $\v$ is not primitive.
\end{Prop}

Then Proposition \ref{prop:tot-wall} and Proposition
\ref{prop:isom} imply that
the following corollary holds.
\begin{Cor}\label{cor:non-primitive}
Assume that $\v$ is not primitive.
Then $\Phi(E)[k]$ is a stable sheaf for a general $E \in {\cal M}_H(\v)$, where $k=0$ or $1$ according as 
$a>0$ or $a \leq 0$.
\end{Cor}

\section{Cohomology Groups of a General Stable Sheaf}\label{sec-cohomology}

In this section, we prove a characterization of those moduli spaces for which weak Brill-Noether fails.  In each of these cases, we compute the cohomology groups of a general stable sheaf in the moduli space.
\subsection{Classifying counterexamples}
We begin by proving a necessary condition for weak Brill-Noether to fail.  Recall from \cref{Def:IsotropicDecomposition} that a Mukai vector $\v=(r,\xi,a)$ admits an isotopric decomposition if $\v= \ell_1 \v_1 + \ell_2 \v_2$, where $\v_i = (r_i, \xi_i, a_i)$ with $$r_i >0, \quad (\xi_i \cdot H)(\xi\cdot H) >0, \quad a_1 > 0, \quad a_2<0,$$ 
$$\v_1^2=\v_2^2 = 0, \quad \langle\v_1, \v_2 \rangle=1, \quad \mbox{and} \quad \{ \ell_1,\ell_2 \}=\left\{ \frac{\v^2}{2},1\right\}.$$
\begin{Thm}\label{thm:wbn}
Let $(X,H)$ be an abelian surface with $H$ a generic polarization. Let  $\v=(r,\xi,a)$ be a Mukai vector such that $\v^2 \geq 0$, and 
$r>0$. 
\begin{enumerate}
    \item If $\v\ne(r,0,-1)e^\eta$ or $\v\ne(1,0,-l)e^\eta$ for some $\eta\in\NS(X)$ and weak Brill-Noether does not hold for $M_{X,H}(\v)$, then $\v$ admits an isotropic decomposition.
    \item If $\v=(r,0,-1)e^\eta$ for $r\ge 2$ or $\v=(1,0,-l)e^\eta$ for some $\eta\in\NS(X)$ and weak Brill-Noether does not hold for $M_{X,H}(\v)$, then $\eta\cdot H<0$, $\eta^2>0$, and in the latter case $\v^2>0$.
\end{enumerate}
\end{Thm}
It follows that in some cases, we can always say that weak Brill-Noether holds.
\begin{Cor}\label{Cor:weak}
Let $\v=(r,\xi,a)$ be a Mukai vector such that $\v^2\geq 0$,
$r>0$ and $\xi \cdot H>0$. 
\begin{enumerate}
       \item If $\v$ is not primitive, then weak Brill-Noether holds for $M_{X, H}(\v)$.
       \item If $\v$ is primitive and $\v^2 + 2 > 2r$, then weak Brill-Noether holds for $M_{X, H}(\v)$.
       \item If $\v_j = (r, \xi, a-j)$ for an integer $j >0$, then weak Brill-Noether holds for $M_{X, H}(\v_j)$.
\end{enumerate}
\end{Cor}
\begin{proof}
(1) By \cref{thm:wbn}(1), if weak Brill-Noether fails for $M_H(\v)$, then $\v$ admits an isotropic decomposition.  Hence, one of the $\v_i$  satisfies $\langle\v,\v_i\rangle=1$, which is impossible if $\v$ is not primitive.  

(2) By \cref{thm:wbn}(1), the existence of an isotropic decomposition implies that 
$$r=\ell_1 r_1+\ell_2 r_2 \geq \ell_1+\ell_2=\ell+1=\frac{\v^2}{2}+1,$$
which cannot happen if $\v^2 + 2 > 2r$.

(3) This follows from part (2) since 
$$\frac{\v_j^2}{2}+1=\frac{\v^2}{2}+rj+1\ge rj+1\ge r+1.$$
\end{proof}
We break the proof of \cref{thm:wbn} into a number of smaller cases covered in \cref{thm:weak,Prop:NonprimitiveCase,Prop:NegativeSlope,Prop:IdealSheaves} below.  In each case, we compute the cohomology of the general sheaf in its moduli space. Since moduli spaces of semi-homogeneous sheaves always satisfy weak Brill-Noether, we may assume that $\v^2 >0$.
\subsubsection{The case when $\v$ is not primitive}
\begin{Prop}\label{Prop:NonprimitiveCase}
Let $\v=(r,\xi,a)$ be a Mukai vector such that $\v^2>0$ and $r>0$.  Then for any $m>1$, $M_H(m\v)$ satisfies weak Brill-Noether.  
\end{Prop}
\begin{proof}
    By \cref{cor:non-primitive}, for the general $E\in M_H(m\v)$, we have $\Phi(E)[k]\in\Coh(\widehat{X})$ for $k=0$ if $a>0$, and $k=1$ if $a\le 0$.  
    By \cref{Lem:FundamentalLemmaAbelian}, we have $H^i(X,E)=0$ for $i>k$.  If $k=0$, then we immediately get that $M_H(m\v)$ satisfies weak Brill-Noether.  If $k=1$, then we get that the general $E\in M_H(m\v)$ has $H^2(X,E)=0$.  Since $m\v$ is not primitive, by \cref{Lem:ReducibleModuli},  $M_{X, H}(m\v)$ is irreducible even when $H$ is not generic. Hence, a general member of $M_{X, H}(m\v)$ remains semistable for a small perturbation of $H$, so we may assume that $H$ is generic and thus that the general $E\in M_H(m\v)$ is a $\mu$-stable vector bundle by \cref{lem:mu-s}.  By Serre duality $h^0(X,E)=h^2(X,E^\vee),$ which vanishes since $E^\vee\in M_H(m(r,-\xi,a))$ is also general and has the same $a$.  It follows that $M_H(m\v)$ satisfies weak Brill-Noether, as required.
\end{proof}
\subsubsection{The case $\xi\cdot H>0$}
\begin{Prop}\label{thm:weak}
Let $\v=(r,\xi,a)$ be a Mukai vector such that $\v^2\geq 0$,
$r>0$ and $\xi \cdot H>0$.  If weak Brill-Noether fails for $M_H(\v)$, then $\v$ admits an isotropic decomposition.  In this case for a general $E\in M_H(\v)$, 
$$h^0(X,E)=\ell_1 a_1,\quad h^1(X,E)=-\ell_2 a_2,\quad\text{ and }\quad h^2(X,E)=0.$$
\end{Prop}

\begin{proof}
We may assume that $\v=(r,\xi,a)$ is primitive by \cref{Prop:NonprimitiveCase}.
We set $\ell:=\v^2/2$ and let $E\in M_H(\v)$ be a general stable sheaf.  
Since $H^2(X,E\otimes\PP_t^\vee)=\Hom(E,\PP_t)^\vee=0$ by the stability of $E$ and $\xi\cdot H>0$, we have $\Phi^2(E)=0$.  If $\Phi^1(E)=0$ for generic such $E$, then $\Phi(E)\in\Coh(\widehat{X})$, so by \cref{Lem:FundamentalLemmaAbelian} $M_H(\v)$ would satisfy weak Brill-Noether.  Thus we have $\Phi^1(E)\ne 0$.  In fact, for the same reason we may even assume that for generic such $E$, $\rk\Phi^1(E)>0$.  Indeed, otherwise $H^i(E \otimes P_y)=0$ for a general $y \in
\widehat{X}$ by cohomology and base change and $M_H(\v)$ would then satisfy weak Brill-Noether, contrary to hypothesis.
By the spectral sequence
$$
E_2^{p,q}=\widehat{\Phi}^p(\Phi^q(E)) \Longrightarrow
\begin{cases}
E & p+q=2\\
0 & p+q \ne 2
\end{cases},
$$
$E$ fits in an exact sequence
\begin{equation}
0 \to \widehat{\Phi}^0 (\Phi^1(E)) \to \widehat{\Phi}^2(\Phi^0(E)) \to E \to \widehat{\Phi}^1(\Phi^1(E)) \to 0.
\end{equation}

Since $\Phi(E)$ is a two term complex of locally free sheaves, 
$\rk \Phi^0(E)=0$ implies $\Phi^0(E)=0$.  In that case, $\Phi(E)[1]\in\Coh(\widehat{X})$, so by \cref{Lem:FundamentalLemmaAbelian} $M_H(\v)$ would satisfy weak Brill-Noether, contrary to hypothesis.  Thus $\rk \Phi^0(E)>0$.  Furthermore, assume that $\rk \Phi^0(E)>0$ and $\rk \Phi^1(E)>0$.
Then Theorem \ref{thm:abel} implies that
$\Phi^0(E)$ and $\Phi^1(E)$ are semistable sheaves with
\begin{equation}
\begin{split}
\v(\Phi^0(E))=& \ell_0 \w_0,\\
\v(\Phi^1(E))=& \ell_1 \w_1,
\end{split}
\end{equation}
where $\w_0^2=\w_1^2=0$, 
$\langle\w_0,\w_1\rangle=-1$, $\{\ell_0,\ell_1\}=\{1,\ell\}$ and $\ell=\v^2/2$.
We set 
$$\widehat{\Phi}(\w_0)=(r_0,\xi_0,a_0),\qquad\widehat{\Phi}(\w_1)=(r_1,\xi_1,a_1),$$
where we note that $a_i=\rk\Phi^i(E)/\ell_i>0$.
Applying \cref{lem:semihom} to $\Phi^1(E)$, we see that $\widehat{\Phi}^0 (\Phi^1(E))=0$ or $\widehat{\Phi}^1 (\Phi^1(E))=0$.
Thus we have two possibilities:
\begin{enumerate}
\item
$\widehat{\Phi}^1 (\Phi^1(E))=0$ and $E$ fits in an exact sequence
$$
0 \to \widehat{\Phi}^0 (\Phi^1(E)) \to \widehat{\Phi}^2(\Phi^0(E)) \to E \to 0.
$$
\item
$\widehat{\Phi}^0 (\Phi^1(E))=0$ and $E$ fits in an exact sequence
$$
0 \to \widehat{\Phi}^2(\Phi^0(E)) \to E \to \widehat{\Phi}^1 (\Phi^1(E))  \to 0.
$$
\end{enumerate}

Assume that $\widehat{\Phi}^1 (\Phi^1(E))=0$.  Then 
\begin{equation*}
\begin{split}
\v(\widehat{\Phi}^0(\Phi^1(E)))=& \ell_1 \u_1,\;\u_1=\widehat{\Phi}(\w_1)=(r_1,\xi_1,a_1),\\
\v(\widehat{\Phi}^2(\Phi^0(E)))=& \ell_0 \u_0,\;\u_0=\widehat{\Phi}(\w_0)=(r_0,\xi_0,a_0).
\end{split}
\end{equation*}
Then $\langle \u_0,\u_1 \rangle=\langle\w_0,\w_1\rangle=-1$.
Since $r>0$, we see that $r_0>0$ and $r_1>0$.
By Lemma \ref{lem:semihom},
$-\xi_1$ is ample and $\xi_0$ is ample.
Hence 
$$1=\langle \u_0,\u_1 \rangle=\xi_0 \cdot \xi_1-r_0 a_1 -r_1 a_0 \leq -3,$$
a contradiction.
Therefore case (1) does not occur.

Now assume that $\widehat{\Phi}^0 (\Phi^1(E))=0$.
Then
\begin{equation*}
\begin{split}
\v(\widehat{\Phi}^2(\Phi^0(E)))=& \ell_0 \u_0,\,\u_0=\widehat{\Phi}(\w_0)=(r_0,\xi_0,a_0),\\
\v(\widehat{\Phi}^1(\Phi^1(E)))=& \ell_1 \u_1,\;\u_1=-\widehat{\Phi}(\w_1)=(-r_1,-\xi_1,-a_1).
\end{split}
\end{equation*}
Since $r>0$, we have $r_0=\rk\widehat{\Phi}^2(\Phi^0(E))>0$.
By Lemma \ref{lem:semihom}, $\xi_0$ is ample.
In particular $\xi_0 \cdot H>0$.
We suppose first that $r_1=0$. Then $-\xi_1$ is an effective divisor with 
$\xi_1^2=2r_1 a_1=0$.
Thus $(-\xi_1)\cdot\xi_0\ge 0$.  
Since $a_1>0$ and $r_0>0$, it follows that
$$1=-\langle\w_0,\w_1\rangle=\langle\u_0,\u_1\rangle=-\xi_0\cdot\xi_1+r_0a_1\ge 1,$$
forcing $\xi_1\cdot\xi_0=0$ and $r_0=a_1=1$.  
But then $\xi_1=0$, so $-a_1>0$, a contradiction.
Therefore $-r_1 > 0$.
Then the stability of $E$ and $\xi \cdot H>0$ imply $(-\xi_1) \cdot H>0$.
Letting $\v_1=\u_0$ and $\v_2=\u_1$ gives the claimed decomposition of $\v$, and taking the long exact sequence on cohomology determines the cohomology groups.
\end{proof}
\subsubsection{The case $\xi \cdot H \leq 0$}\label{sec:NegativeSlope}

\begin{Prop}\label{Prop:NegativeSlope}
Let $\v=(r,\xi,a)$ be a Mukai vector such that $\v^2\geq 0$,
$r \geq 2$ and $\xi \cdot H\leq 0$, and let $H$ be a generic polarization with respect to $\v$.
\begin{enumerate}
\item
Assume that $\v \ne (r,0,-1)e^{\eta}$ $(\eta \in \NS(X))$.
If $M_H(\v)$ does not satisfy the weak Brill-Noether property, then $\v$ admits an isotropic decomposition.  In this case for a general $E\in M_H(\v)$, 
$$h^0(X,E)=0,\quad h^1(X,E)=-\ell_1 a_1,\quad\text{ and }\quad h^2(X,E)=\ell_2 a_2.$$
\item
Assume that $r \mid \xi$ and $\v =(r,0,-1)e^{\xi/r}$.
If $M_H(\v)$ does not satisfy the weak Brill-Noether property, then
$\xi^2>0$ and $\xi\cdot H<0$.   In this case for a general $E\in M_H(\v)$,
$$h^0(X,E)=0,\quad h^1(X,E)=1,\quad\text{ and }\quad h^2(X,E)=a+1.$$
\end{enumerate}
\end{Prop}
\begin{proof}
(1)  Since $\v \ne (r, 0, -1)e^{\eta}$ and $H$ is generic, the general $E\in M_{X,H}(\v)$ is a $\mu_H$-stable vector bundle by Lemma \ref{lem:mu-s}.
If $\xi \cdot H =0$, then by the Hodge Index Theorem, $\xi^2 \leq 0$. Hence, $\v^2 = \xi^2-2ra \geq 0$ implies that $\chi(\v)=a \leq 0$. 
By $\mu_H$-stability $h^0(E)=0$ and by Serre duality and $\mu_H$-stability $h^2(E)=0$. We conclude that $h^1(E)=-a$ and weak Brill-Noether holds, contrary to our assumption.  Thus we must have $\xi \cdot H <0$, from which it follows that  $h^0(E)=0$ by stability.
By Serre duality $h^1(E)= h^1(E^{\vee})$ and $h^2(E)=h^0(E^\vee)$. The Mukai vector $\v^{\vee}= (r, - \xi, a)$ of $E^{\vee}$ satisfies $\xi \cdot H >0$. If weak Brill-Noether fails for $M_H(\v)$, then it also fails for $M_H(\v^\vee)$, so the decomposition from \cref{thm:weak} for $\v^{\vee} = \ell_1 \v_1 + \ell_2 \v_2$ with $\v_i = (r_i, \xi_i, a_i)$ yields a decomposition for $\v= \ell_1 \v_1' + \ell_2 \v_2'$ with $\v_i ' = (r_i, -\xi_i, a_i)$ satisfying the conditions of the theorem.

(2) Now assume that $r \mid \xi$ and $\v =(r,0,-1)e^{\xi/r}$.
We set $\eta:=\xi/r$.
Then a general $E \in M_H(\v)$ fits in an exact sequence 
$$ 
0 \to E \to \oplus_{j=1}^r L_j \to k(x) \to 0,
$$
where $L_j$ are general line bundles with
$c_1(L_j)=\eta$ and $k(x)$ is the structure sheaf of a point $x \in X$
(see \cite[Thm. 2.20, Prop. 2.21, Cor. 4.5]{Mu:Fourier}).
Suppose first that $\xi\cdot H=0$ so that $\eta\cdot H=0$.  Then $h^0(L_j)=h^2(L_j)=0$ for all $j$ since $H$ is ample and the $L_j$ are general.  Hence, $h^0(E)=h^2(E)=0$, so $M_H(\v)$ satisfies weak Brill-Noether, contrary to assumption.  Thus we must have $\eta\cdot H<0$ so that $h^0(E)=h^0(L_j)=0$ for all $j$ by stability.  Since $h^2(E)=h^2(\oplus_{j=1}^r L_j)$, and we are assuming that weak Brill-Noether fails, we have $\eta^2>0$ so $\xi^2>0$ as required.
\end{proof}
\begin{Prop}\label{Prop:IdealSheaves}
    Let $\v=(1,\xi,a)$  be a Mukai vector such that $\v^2\ge 0$.  Assume that $\xi\cdot H\le 0$.  If $M_H(\v)$ fails to satisfy weak Brill-Noether, then $\xi^2>0$, $\xi\cdot H<0$, and $\v^2>0$.  In this case, for general $E\in M_H(\v)$,
    $$h^0(X,E)=0,\quad h^1(X,E)=\frac{\xi^2}{2}-a,\quad\text{ and }\quad h^2(X,E)=\frac{\xi^2}{2}$$
\end{Prop}
\begin{proof}
    Write $\v=(1,0,-l)e^\xi$ for $l\ge 0$, and let $L$ be a general line bundle with $c_1(L)=\xi$.  Supposing $\xi^2\le 0$, we have $h^0(L)=h^2(L)=0$ so the same holds true for $E=I_Z\otimes L$ the general element of $M_H(\v)$.  But then $M_H(\v)$ satisfies weak Brill-Noether, contrary to hypothesis.  Thus $\xi^2>0$, so in particular  $\xi\cdot H<0$.  It follows that $h^0(L)=0=h^1(L)$, while $h^2(L)=\frac{\xi^2}{2}=a+l$.  Thus $E=I_Z\otimes L$ satisfies $h^1(E)=l$ and $h^2(E)=a+l$.  Since $M_H(\v)$ fails to satisfy weak Brill-Noether, we must have $l>0$.  This is equivalent to $\v^2>0$, as claimed.
\end{proof}

\subsection{Constructing Examples}\label{sec:ConstructingExmples}
By \cref{thm:wbn}, for $\v=(r,\xi,a)$ with $r>0$ and $\v^2\ge 0$ either the general sheaf $E\in M_H(\v)$ has at most one nonzero cohomology group determined as in \cref{Rem:WBNDetermingsCohomology} or it falls into one of three cases:
\begin{enumerate}
    \item $\v$ admits an isotropic decomposition.
\item $\v=(r,0,-1)e^\eta$ with $r\ge 2$ and $\eta\in\NS(X)$ such that $\eta\cdot H<0$ and $\eta^2>0$.
\item $\v=(1,0,-l)e^\eta$ with $l>0$ and $\eta\in\NS(X)$ such that $\eta\cdot H<0$ and $\eta^2>0$.
\end{enumerate}

In the latter two cases, it is clear that such counterexamples exist and how to construct them.  So we consider the first case of a Mukai vector admiting an isotropic decomposition and construct counterexamples in the following proposition.
\begin{Prop}\label{Prop:ConstructingCounters}
Assume that there is a pair of isotropic Mukai vectors $\v_1=(r_1,\xi_1,a_1)$ and $\v_2=(r_2,\xi_2,a_2)$ satisfying 
\begin{equation}\label{eq:pairs}
r_1, r_2>0,\quad a_1>0,\quad a_2<0,\quad \langle \v_1,\v_2 \rangle=1.
\end{equation}
Then there exists a polarization $H$ such that $M_H(\ell_1 \v_1+\ell_2 \v_2)$ with 
$\{\ell_1,\ell_2 \}=\{\ell,1 \}$ fails weak Brill-Noether.
\end{Prop}
\begin{proof}
Since
$(r_2 \xi_1-r_1 \xi_2)^2=-2r_1 r_2 \langle \v_1,\v_2 \rangle<0$,
there is an ample divisor $L$ with $L \cdot (r_2 \xi_1-r_1 \xi_2)=0$.
Replacing $\xi_i$ by $-\xi_i$, we assume that $\xi_i \cdot L>0$.
We take a ${\Bbb Q}$-divisor $\epsilon$ such that
$H=L-\epsilon$ is ample and
$H \cdot (r_2 \xi_1-r_1 \xi_2)<0$.
For $E \in M_H(\ell_1 \v_1+\ell_2 \v_2)$ with 
$\{\ell_1,\ell_2 \}=\{\ell,1 \}$, we have an exact sequence
$$
0 \to E_1 \to E \to E_2 \to 0 
$$
where $E_i \in  M_H(\ell_i \v_i)$.
Hence $h^0(E) h^1(E) \ne 0$ for all $E \in M_H(\ell_1 \v_1+\ell_2 \v_2)$.
Thus as along as we can choose a pair of isotropic Mukai vectors $\v_1$ and $\v_2$ as above, we can construct counterexamples to weak Brill-Noether.    
\end{proof}

In the next section, we rephrase the existence of this pair entirely in terms of the Neron-Severi group of $X$, and characterize precisely when counterexamples exist.

\section{Reinterpreting counterexamples in terms of $\NS(X)$}\label{sec-Neron-Severi}
In this section we prove the main thoerem \cref{Thm:MainTheorem1}.  We first reformulate the criterion in \cref{thm:wbn} giving counterexamples to weak Brill-Noether in terms of the Neron-Severi group.
\begin{Lem}\label{Lem:Translation}
The existence of a pair of isotropic Mukai vectors $\v_1=(r_1,\xi_1,a_1)$ and $\v_2=(r_2,\xi_2,a_2)$ satisfying 
$$r_1,r_2>0, a_1<0, a_2>0, \langle \v_1,\v_2 \rangle=1$$
is equivalent to the existence of $A, B \in \NS(X)$ and positive integers $r_1, r_2$ such that 
$$A^2 > 0,\quad B^2<0,\quad 2r_2\mid A^2,\quad 2r_1\mid B^2, \quad \text{and} \quad (r_1 A - r_2 B)^2 = -2r_1r_2.$$
\end{Lem}
\begin{proof}
If $\v_1$ and $\v_2$ are isotropic and satisfy the above conditions, then let the $r_i$'s be the same and set $A=\xi_2$ and $B=\xi_1$.  Then 
$$(r_1A-r_2B)^2=(r_1\xi_2-r_2\xi_1)^2=-2r_1r_2\langle\v_1,\v_2\rangle=-2r_1r_2.$$
Furthermore, $A^2=\xi_2^2=2a_2r_2>0$ and $B^2=2a_1r_1<0$.  

Conversely, if $r_1$, $r_2$ and $A,B$ are  as above, then set $a_2=\frac{A^2}{2r_2}>0$ and $a_1=\frac{B^2}{2r_1}<0$.  We then have that  $\v_1=(r_1,B,a_1)$ and $\v_2=(r_2,A,a_2)$ are isotropic Mukai vectors satisfying the required conditions.
\end{proof}
In particular this immediately gives the following result in Picard rank one.
\begin{Cor}\label{Cor:WBNPicardRankOne}
Assume that $X$ is an abelian surface with $\NS(X)=\Z H$.  Let $\v=(r,dH,a)$ be a Mukai vector such that $r>0$, $d>0$, and $\v^2\geq 0$. Then $M_H(\v)$ satisfies weak Brill-Noether.
\end{Cor}
\begin{proof}
Since $d>0$, by \cref{thm:wbn} the only way $M_H(\v)$ could fail weak Brill-Noether is if $\v$ admits an isotropic decomposition as in \cref{thm:wbn}, but this is equivalent to the existence of $A,B\in\NS(X)$ and $r_1,r_2\in\N$ as in \cref{Lem:Translation}.  Since $\rho(X)=1$, the intersection pairing on $\NS(X)$ is positive definite so that $(r_1 A-r_2 B)^2\ge 0$.
\end{proof}
We proceed by considering the higher Picard rank case.  We begin with a couple of purely arithmetic lemmas.
\begin{Lem}\label{Lem:ArithmeticTrick}
Let $r_1,r_2,x_1,x_2,y_1, y_2$ be integers such that $r_1,r_2>0$ and 
\begin{equation}
r_1 r_2
= -(r_1 x_2-r_2 x_1)(r_1 y_2-r_2 y_1).
\end{equation}
Then $x_1 x_2 \geq 0$ and $y_1 y_2 \geq 0$.
\end{Lem}
\begin{proof}
We set $k:=\gcd(r_1,r_2)>0$.
Then 
$$
r_1=k r_1', r_2=k r_2', \gcd(r_1',r_2')=1
$$
and
\begin{equation}
(r_1' x_2-r_2' x_1)(r_1' y_2-r_2' y_1)=-r_1' r_2'.
\end{equation}
By $\gcd(r_1',r_2')=1$, 
there are positive integers $s_1,s_2,t_1,t_2$ such that
\begin{equation}
\begin{split}
(r_1' x_2-r_2' x_1)=& \pm s_1 s_2,\; s_1 \mid r_1', s_2 \mid r_2'\\
(r_1' y_2-r_2' y_1)=& \mp t_1 t_2,\; t_1 \mid r_1', t_2 \mid r_2'.
\end{split}
\end{equation}
Then we see that
$$
\frac{x_1}{s_1}, \frac{x_2}{s_2}, \frac{y_1}{t_1},\frac{y_2}{t_2} \in {\Bbb Z}
$$ 
and
\begin{equation}
\begin{split}
\frac{r_1'}{s_1}\frac{x_2}{s_2}-\frac{r_2'}{s_2}\frac{x_1}{s_1}=& \pm 1 \\
\frac{r_1'}{t_1}\frac{y_2}{t_2}-\frac{r_2'}{t_2}\frac{y_1}{t_1}=& \mp 1.
\end{split}
\end{equation}
Hence $x_1 x_2 \geq 0$ and $y_1 y_2 \geq 0$.
\end{proof}
In the following we denote by $U$ the standard hyperbolic lattice.
\begin{Lem}\label{Lem:CriterionSubHyperbolic}
Let $L$ be an even lattice of signature $(1,1)$ and discriminant $d(L)$. Then the following conditions are equivalent.
\begin{enumerate}
\item[(i)]
$\sqrt{-d(L)} \in {\Bbb Q}$.
\item[(ii)]
$L$ is isomorphic to a sublattice of $U$.
\item[(iii)]
$L$ contains an isotropic vector.
\end{enumerate}  
\end{Lem}

\begin{proof}
Let 
$$
\begin{pmatrix}
2a & b\\
b & 2c
\end{pmatrix}, (a,b,c \in {\Bbb Z})
$$
be the intersection matrix of $L$.

(i) $\implies$ (ii).
Assume that $d(L)=4ac-b^2=-n^2$, $n \in {\Bbb Z}$.
Since $\frac{b-n}{2}, \frac{b+n}{2} \in {\Bbb Z}$ and
$ac=\frac{b-n}{2} \frac{b+n}{2}$,
there are relatively prime integers
$x,y$ such that 
\begin{equation}
(a,\tfrac{b+n}{2})=(sx,sy),\; 
(\tfrac{b-n}{2},c)=(tx,ty),
\; (s, t \in {\Bbb Z}).
\end{equation}
Since $b=sy+tx$, we get
\begin{equation}\label{eq:L}
\begin{pmatrix}
2a & b\\
b & 2c
\end{pmatrix}
=
\begin{pmatrix}
s & x\\
t & y
\end{pmatrix}
\begin{pmatrix}
0 & 1\\
1 & 0
\end{pmatrix}
\begin{pmatrix}
s & t\\
x & y
\end{pmatrix}.
\end{equation}
Hence $L$ is isomorphic to a sublattice of $U$.


(ii) $\implies$ (iii).
We may assume that $L$ is a sublattice of $U$.
Then $d(L)U \subset L$. Hence there is an isotropic vector in $L$.

(iii) $\implies$ (i).
By our assumption, we may assume that $c=0$. Then we have $d(L)=-b^2$.
Therefore our claim holds.
\end{proof}
We study the existence of $(A,B,r_1,r_2)$ as in \cref{Lem:Translation} by reducing to the case of rank two sublattices of $\NS(X)$ that contain an ample class.  Let $L= \Z H + \Z D\subset\NS(X)$ be a rank two sublattice containing an ample class $H$ with intersection pairing
$$H^2 = 2n, \quad D^2 = 2m, \quad \mbox{and} \quad H \cdot D = k.$$
We would like to find classes $A, B \in L$ and positive integers $r_1, r_2$ such that $A^2 > 0$, $B^2<0$, $2r_2\mid A^2$, $2r_1\mid B^2$, and $(r_1 A - r_2 B)^2 = -2r_1r_2$.
\begin{Prop}\label{Thm:WBNPicardRankTwo}
Let $L$ be as above and write $\Delta=-d(L)$.
\begin{enumerate}
    \item If $\Delta$ is a perfect square, then there does not exist such $(A, B, r_1, r_2)$.
    \item If $\Delta$ is not a perfect square, then there exists infinitely many such $(A, B, r_1, r_2)$.
\end{enumerate}
\end{Prop}
\begin{proof}
Express $A= aH + bD$ and $B= cH + dD$ for integers $a,b,c,d$.
    Then the (in)equalities become 
    \begin{eqnarray*}
    A^2 &=& 2na^2 +2mb^2+ 2k ab >0 \\
    B^2 &=& 2nc^2 + 2md^2 + 2k cd <0 \\
    (r_1A - r_2B)^2 &=& r_1^2 (2n a^2 + 2mb^2 + 2k ab) + r_2^2(2nc^2 + 2md^2 + 2k cd) \\ & \ & - 2r_1r_2(2nac + 2mbd + k(ad+bc))= - 2r_1 r_2.
    \end{eqnarray*}
Making the substitutions $x= r_1 a - r_2 c$ and $y = r_1 b -r_2 d$, the last equation can be written as 
$$nx^2 + k xy + m y^2 = -r_1r_2.$$
By the Hodge Index Theorem,  the signature of the quadratic form is $(1, 1)$, so $\Delta:=k^2 - 4mn >0$. 

Assume first that $\Delta$ is a square.  Then by \cref{Lem:CriterionSubHyperbolic} it follows that we have an isometric embedding $\psi\colon L\to U=\Z f+\Z g$.  Assume to the contrary that there exists a tuple $(A,B,r_1,r_2)$ as above.  Set 
$$\psi(A)=x_2f+y_2g$$
$$\psi(B)=x_1f+y_1g$$
Then 
$$r_1r_2=-\frac{1}{2}\prs{r_1A-r_2B}^2=-(r_1x_2-r_2x_1)(r_1y_2-r_2y_1),$$
so by \cref{Lem:ArithmeticTrick} it follows that $x_1x_2\ge 0$ and $y_1y_2\ge 0$, and thus $x_1x_2y_1y_2\ge 0$.  On the other hand, because $\psi(A)^2>0$ and $\psi(B)^2<0$, we have
$$x_1y_1x_2y_2=\frac{\psi(A)^2}{2}\frac{\psi(B)^2}{2}<0,$$
a contradiction.

Now assume that $\Delta$ is not a square.  We would like to show that there exist solutions to this system.  First, observe that the quadratic form represents negative numbers. This is clear if $n=m=0$ or one of $n$ or $m$ is negative. Hence, with out loss of generality, we may assume that $n >0$. The value of the quadratic form at $x= -k$, $y= 2n$ is $n(-k^2 + 4mn)$, which is negative by assumption. Consequently, the quadratic form represents infinitely many negative numbers. Fix  a negative number $-R$ that the form represents. We set $r_1=1$ and $r_2=R$.

We can generate more solutions of these equations using the Pell's equation.
Multiply the equation by $4n$ to obtain
$$4n^2 x^2 + 4nk xy  + 4nm y^2 = -4Rn.$$
Rewrite the equation as follows
$$4n^2 x^2 + 4nk xy + k^2 y^2 + (4mn - k^2)y^2= - 4Rn.$$
Making the substitution $u = 2n x + ky$ and $v= y$, we get the equation
$$u^2 - \Delta v^2 = -4Rn,$$ where $\Delta= k^2 - 4mn$ is the discriminant. Suppose we have a solution of this equation $(u_0, v_0)$.

Given a solution $(u_1, v_1)$ to the Pell's equation
$$u^2 - \Delta v^2 =1,$$ we get a new solution of the equation 
$$u^2 - \Delta v^2 = -4Rn$$ of the form
$$(u_0u_1 + \Delta v_0 v_1, u_0 v_1 + u_1 v_0).$$

Since $\Delta$ is not a square, the Pell's equation has infinitely many solutions, in particular solutions where the absolute value tends to infinity. 

Now let $(x_0,y_0)$ and $(c,d)$ be solutions of 
$nx^2 + k xy + m y^2= -R$ such that 
$$d\sim \frac{-k-\sqrt{\Delta}}{2m}c,\quad y_0\sim\frac{-k+\sqrt{\Delta}}{2m}x_0,$$
where $|c|,|x_0|\gg 0$.  We will choose their signs appropriately in a moment.
Letting  $a= x_0 + Rc$ and $b = y_0 + Rd$, we have that
\begin{align*}
na^2+kab+mb^2&=    n(x_0+Rc)^2 + k (x_0+Rc)(y_0+Rd) + m (y_0+Rd)^2= nx_0^2 + k x_0y_0 + my_0^2 \\
&+ nR^2 c^2 + kR^2 cd + mR^2 d^2 + 2nx_0Rc + kRy_0c + k R x_0d + 2my_0Rd\\
&=-R-R^3+2nx_0Rc + kRy_0c + k R x_0d + 2my_0Rd\\
&\sim -R-R^3+2nx_0cR+kcx_0 R\prs{\frac{-k+\sqrt{\Delta}}{2m}}+kcx_0 R\prs{\frac{-k-\sqrt{\Delta}}{2m}}\\
&+2mcx_0 R\prs{\frac{-k+\sqrt{\Delta}}{2m}\frac{-k-\sqrt{\Delta}}{2m}}\\
&=-R-R^3+\prs{2n+\frac{-k^2+k\sqrt{\Delta}}{2m}+\frac{-k^2-k\sqrt{\Delta}}{2m}+\frac{k^2-\Delta}{2m}}cx_0 R\\
&=-R-R^3+\prs{\frac{4mn-k^2-\Delta}{2m}}cx_0 R=-R-R^3-\frac{\Delta}{m}cx_0 R.
\end{align*}
Thus depending on the sign of $m$, we can choose $c$ and $x_0$ arbitrarily large and of the same sign (if $m<0$) or arbitrarily large and of opposite signs (if $m>0$) to make $na^2+kab+mb^2>0$. 
For any choice of the infinitely many $(x_0,y_0)$ and $(c,d)$ satisfying these conditions, we get a $4$-tuple $(A,B,r_1,r_2)=\prs{\prs{x_0+Rc}H+\prs{y_0+Rd}D,cH+dD,1,R}$ of the type required.  Observe that $r_1=1$ and $r_2=R$, so the divisibility conditions are automatically satisfied by construction.
\end{proof}
This result is already enough to classify when counterexamples to weak Brill-Noether exist in Picard rank two.
\begin{Cor}\label{Cor:WBNPicardRankTwo}
Let $X$ be an abelian surface with $\rho(X)=2$.  If $-d(\NS(X))$ is a square, then weak Brill-Noether holds for all moduli spaces $M_H(\v)$ with $\v=(r,\xi,a)$ such that $r>0$, $\v^2\ge 0$ and polarization $H$ such that $\xi\cdot H>0$ .  Otherwise, for infinitely many $\v=(r,\xi,a)$ such that $r>0$, $\v^2\ge 0$, there exists a polarization $H$ with $\xi\cdot H>0$ such that $M_H(\v)$ does not satisfy weak Brill-Noether.
\end{Cor}
In even higher Picard rank, we begin with a guiding example on how we find a rank two sublattice giving rise to $(A,B,r_1,r_2)$ as in \cref{Lem:Translation}.
\begin{Ex}\label{Ex:ProductEllipticCurves}
    Let $C$ be an elliptic curve and set $X:=C \times C$.
Let $f$ and $g$ be fibers of the two projections $X \to C$ and
$\Delta$ the diagonal of $X=C \times C$. 
Then 
$$
\NS(X)\supset{\Bbb Z}f+{\Bbb Z}g+{\Bbb Z}D, 
$$
where $D=\Delta-f-g$.
We take an ample divisor $H:=f+mg$ with $\sqrt{m} \not \in {\Bbb Z}$.
Let $(s,t)$ be a pair of positive integers such that $t^2-ms^2=1$ and $s>1$.
We set $\eta:=-sH+tD$, which satisfies $(\eta^2)=-2$.
Then setting $A=H$, $B=\eta+H$ and $r_1=r_2=1$, we have
$$A^2=H^2=2m>0, B^2=(\eta+H)^2=-2((2s-1)m+1)<0, (r_1A-r_2B)^2=(-\eta)^2=-2=-2r_1r_2.$$
Moreover, the same is true for all of the infinitely many pairs of positive integers $(s,t)$ such that $t^2-ms^2=1$ with $s>1$.
\end{Ex}
\begin{Prop}\label{Prop:InfinitelyManyCounterexamples}
    Let $X$ be an abelian surface with $\rho(X)\ge 3$.  The subset of $\NS(X)\times\NS(X)\times\N\times\N$ given by
    $$\BB=\Set{(A,B,r_1,r_2)\ |\ A^2>0,B^2<0,2r_2\mid A^2, 2r_1\mid B^2,(r_1A-r_2B)^2=-2r_1r_2}$$
    is infinite.
\end{Prop}
\begin{proof}
Suppose first that $\rho(X)=4$.  Then by \cite[Theorem 4.1]{SM:SingularAbelianSurfaces}, $X$ is isogenous to a self product $C\times C$ of an elliptic curve $C$ with complex multiplication.  On $C\times C$, we have two divisor classes $H$ and $D$ as in \cref{Ex:ProductEllipticCurves}.  Up to multiplication by a positive integer $n$, we may assume that $H$ and $D$ are in $\NS(X)$.  If $n=1$, then we obtain infinitely many elements in $\BB$ as in \cref{Ex:ProductEllipticCurves}.  If $n>1$, then we observe that $-d(\Z nH+\Z nD)=-n^2d(\Z H+\Z D)=4n^2m$ cannot be a perfect square since $m$ is not one.  Then it follows from \cref{Thm:WBNPicardRankTwo} that we again get infinitely many elements in $\BB$.

Thus we may assume that $\rho(X)=3$.  Pick an integral ample divisor class $H$. 
Since $\rho(X)=3$, the nef cone is the round cone 
$$\Nef(X)=\Set{z\in\NS(X)_\R\ |\ z^2\ge 0, z\cdot H>0}.$$
Consider a rational hyperplane $W$ through the origin and containing the ray spanned by $H$.  Then $W\cap\Nef(X)$ is spanned by two isotropic rays that are either both irrational or both rational.  If they are both irrational, then by \cref{Thm:WBNPicardRankTwo} we already get infinitely many elements in $\BB$.  So we may assume that $W\cap\Nef(X)$ is spanned by two rational isotropic rays.  By scaling appropriately we conclude that there exists at least two primitive isotropic divisor classes in $\NS(X)\cap\Nef(X)$ that are linearly independent  over $\R$.  By \cite[Lemma 3.3]{Ma:DecompositionsAbelianSurface}, it follows that $X$ contains an elliptic curve and is thus isogenous to a decomposable abelian surface $X'$ that also has $\rho(X')=3$ since the Picard number is stable under isogenies.  

By \cite[Proposition 4.1]{Ma:DecompositionsAbelianSurface}, we have $\NS(X')\cong U\oplus\langle-2N\rangle$ for some $N\in\N$ and we can find primitive isotropic vectors $x$, $y$, and $z$ in $\NS(X')\cap\Nef(X')$ such that 
$$x\cdot y=1,\quad x\cdot z=1,\quad y\cdot z=N.$$
We claim that the rank two sublattice $L\subset\NS(X')$ spanned by $x+y$ and $x+z$ satisfies $-d(L)$ is not a square. 
Observe first that 
$$(x+y)^2=2,\quad (x+z)^2=2,\quad (x+y)\cdot(x+z)= 2+N.$$
Suppose to the contrary that $-d(L)$ is a square.
Then there exists $d\in\Z$ such that
$$d^2=-d(L)=( 2+N)^2-4=N^2+4N,$$
which is easily seen to be impossible.
Up to multiplication by some $k\in\N$, we may assume that $kL\subset\NS(X)$.  But then $-d(kL)=-k^2d(L)$ cannot be a perfect square since $-d(L)$ is not a perfect square.  It then follows that $\BB$ contains infinitely many elements by \cref{Thm:WBNPicardRankTwo}.
\end{proof}
\begin{proof}[Proof of \cref{Thm:MainTheorem1}]
The Picard rank one case follows from \cref{Cor:WBNPicardRankOne} and the Picard rank two case follows from \cref{Cor:WBNPicardRankTwo}.  The higher Picard rank case follows from \cref{thm:wbn,Prop:InfinitelyManyCounterexamples} and \cref{Prop:ConstructingCounters}.
\end{proof}
\subsection{Infinitely many walls}
Let $X$ be an abelian surface and
$\v=(r,\xi,a)$ a Mukai vector such that $r>0$ and
$\v^2=2\ell>0$.
In proving \cref{Thm:MainTheorem1}, we have shown that if $X$ does not contain an elliptic curve  and $\rho(X)=2$ or when $\rho(X)\ge 3$, then we can always construct infinitely many counterexamples $\v$ that do not possess the weak Brill-Noether property.  In this subsection, we shall refine this idea and prove that if $\v$  admits one isotropic decomposition then it actually admits infinitely many.

\begin{Prop}\label{prop:infinite}
Let $X$ be an abelian surface and $\v=(r,\xi,a)$ be a Mukai vector with 
$\v^2=2\ell>0$. 
Assume that $\v$ admits an isotropic decomposition: there are isotropic Mukai vectors 
$\v_i=(r_i,\xi_i,a_i)$ ($i=1,2$) and a decomposition
$\v=\ell \v_1+\v_2$ such that
$$r_1,r_2>0,\quad a_1>0,\quad a_2<0,\quad \langle\v_1,\v_2\rangle=1.$$
Then:
\begin{enumerate}
\item
There are infinitely many totally semistable walls for Gieseker stability (corresponding to infinitely many isotropic decompositions).
 \item
There are infinitely many chambers of the ample cone containing ample divisors $H$ such that $M_H(\v)$ fails weak Brill-Noether.
\end{enumerate} 
\end{Prop}
The result follows from the following lemmas.

\begin{Lem}
Let $\v_1:=(r_1,\xi_1,a_1)$ be a Mukai vector such that
$0<r_1<r$ and $(r\xi_1-r_1 \xi)^2 < 0$.
Then
$\v_1$ defines a wall for Gieseker stability if and only if
$ \v_1^2\ge 0$ and $(\v-\v_1)^2  \geq 0$.
\end{Lem}

\begin{proof}
We take an ample divisor $H$ such that $H \cdot(r_1 \xi-r_1 \xi)=0$.
Since $\v_1^2\ge 0$ and $(\v-\v_1)^2 \geq 0$,
there are semi-stable sheaves $E_1 \in { M}_H(\v_1)$ and $E_2 \in { M}_H(\v-\v_1)$.
Then $E:=E_1 \oplus E_2$ is a properly $\mu$-semistable sheaf with $\v(E)=\v$. 
\end{proof}
\begin{Rem}
Under the assumptions of \cref{prop:infinite} we must either have (i) $X$ does not contain an elliptic curve and
$\rho(X)=2$ or (ii) $\rho(X) \geq 3$ by \cref{Thm:MainTheorem1}.  We show now that the infinitude of counterexamples is due to the presence of a large automorphism group in these cases.
\end{Rem}

\begin{Lem}\label{lem:O}
Let $L$ be an even lattice of signature $(1,1)$ and
assume that $\sqrt{-d(L)} \not \in {\Bbb Z}$.
Then the isometry group $O(L)$ of $L$ is infinite.
\end{Lem}

\begin{proof}
We take an integral basis $e_1,e_2$ of $L$ and let
$
Q:=
\begin{pmatrix}
2a & b\\
b & 2c
\end{pmatrix}
$
be the intersection matrix of $L$ with respect to this basis.
We set 
$
A:=
\begin{pmatrix}
-b & -2c\\
2a & b
\end{pmatrix}.
$
Then $A^2=(b^2-4ac)I_2$, where $I_2$ is the $2\times 2$ identity matrix, and
$$
{}^t A Q+QA=0,\; {}^t A QA=-(b^2-4ac)Q.
$$
For an integral solution $(t,u)$ of the Pell's equation
$t^2-(b^2-4ac)u^2=4$, we see that
$$
g:=\tfrac{1}{2}(t I_2+u A) \in SL(2,{\Bbb Z})
$$
and $g$
defines an isometry of $L$.
Since this Pell's equation has infinitely many integral solutions, $O(L)$ is infinite.
\end{proof}
From now on,
$\v_1:=(r_1,\xi_1,a_1)$ is a Mukai vector such that
$0<r_1<r$ and $(r\xi_1-r_1 \xi)^2 < 0$.
\begin{Lem}\label{lem:infinite}
Assume that there is a divisor $D$ and a hyperbolic sublattice $L$ of $\NS(X)$
such that
\begin{enumerate}
\item
$\xi+rD, r \xi_1-r_1 \xi \in L$.
\item
$L$ does not contain an isotropic element, that is $\sqrt{-d(L)} \not \in {\Bbb Q}$.
\end{enumerate}
Then 
\begin{enumerate}
\item[(i)]
Assume that $\v_1$ defines a wall for Gieseker stability. Then
there are infinitely many walls for Gieseker stability.
 \item[(ii)]
Assume that $\v_1$ defines a totally semistable wall, that is, $\v_1$ is an isotropic vector that is part of an isotropic decomposition for $\v$, so in particular $\langle\v,\v_1 \rangle=1$. 
Then
there are infinitely many totally semistable walls for Gieseker stability.
Moreover if $\xi \not \in L^\perp$, then each one of these walls bounds a chamber of the ample cone containing ample divisors $H$ such that ${ M}_H(\v)$ fails weak Brill-Noether.
\end{enumerate} 
\end{Lem}

\begin{proof}
(i)
We set 
$$
\xi':=\xi+rD, \quad \xi_1':=\xi_1+r_1 D, \quad \v':=\v e^D,\quad \v_1':=\v_1e^D.
$$
Then $r \xi_1'-r_1 \xi'=r \xi_1-r_1 \xi \in L$.

By our assumption on $L$, it follows from \cref{lem:O} that
there is an infinite order automorphism $\varphi$ of $L$.
Replacing $\varphi$ by $\varphi^N$, we may assume that
$\varphi$ acts trivially on $L/r L$.
Then $\varphi^n(r\xi_1-r_1 \xi) \ne r\xi_1-r_1 \xi $ for $n \ne 0$.
We define $\eta_n \in \NS(X)$ by
$$
\eta_n+rD=\varphi^n(\xi_1')-r_1 \frac{\varphi^n(\xi')-\xi'}{r}.
$$
Then we see that
$$
\eta_n=\frac{1}{r}\varphi^n(r\xi_1-r_1 \xi)+\frac{r_1}{r}\xi.
$$
We note that $\eta_n =\eta_m$ if and only if $n=m$.
We set 
$$
\w_n:=(r_1,\eta_n, b_n)=(r_1,\varphi^n (\xi_1'),a_1)e^{\frac{(\xi'-\varphi^n(\xi'))}{r}-D} \in
H^{2*}(X,{\Bbb Z}).
$$
Since $\varphi^n(\v')e^{\frac{(\xi'-\varphi^n(\xi'))}{r}}=\v'$, we see that 
$$
\langle \w_n,\v \rangle=\langle \v_1,\v \rangle,\;
\langle \w_n,\w_n \rangle=\langle \v_1,\v_1 \rangle.
$$
Hence, if $\v_1$ defines a wall for Gieseker stability, then the $\w_n$ also define walls for Gieseker stability  so there are infinitely many walls for Gieseker stability.

(ii)
Assume that $\v_1$ is an isotropic Mukai vector with $\langle \v_1,\v \rangle=1$ that forms part of an isotropic decomposition: $\v=\ell \v_1+\v_2$.  Then $\v_2$ is an isotropic Mukai vector with 
$\langle \v_1,\v_2 \rangle=1$ and $r_2:=r-\ell r_1>0$.
We set $\v_2=(r_2,\xi_2,a_2)$ and
$\zeta_n:=\xi-\ell \eta_n$.
Then
$$
\zeta_n=-\frac{\ell}{r}\varphi^n(r \xi_1-r_1 \xi)+\frac{r_2}{r}\xi.
$$

We take $e_1, e_2 \in L \otimes {\Bbb R}$ 
such that $(e_1^2)=(e_2^2)=0$ and $(e_1 \cdot e_2)=1$. 
Since $\varphi$ is an isometry of $L$, there is a real number $\alpha \ne 0$ such that 
$$(\varphi(e_1),\varphi(e_2))=(\alpha e_2,\alpha^{-1} e_1) \quad \mbox{or} \quad  
(\varphi(e_1),\varphi(e_2))=(\alpha e_1,\alpha^{-1} e_2).$$
If $(\varphi(e_1),\varphi(e_2))=(\alpha e_2,\alpha^{-1} e_1)$, then $\varphi^2$ is the identity.
Hence $(\varphi(e_1),\varphi(e_2))=(\alpha e_1,\alpha^{-1} e_2)$.
Replacing $\varphi$ by $\varphi^2$, we assume that $\alpha>0$.
We may also assume that $\alpha>1$.
We write 
\begin{equation}
\begin{split}
r \xi_1-r_1 \xi= & xe_1+ye_2,\\
\xi= & ze_1+w e_2+\nu, \nu \in L^\perp \otimes {\Bbb Q}
\end{split}
\end{equation}
By our assumption that $\xi\notin L^\perp$, $z \ne 0$ or $w \ne 0$. 
Then we see that
\begin{equation}
\begin{split}
r\eta_n &=(\alpha^n x+r_1 z)e_1+(\alpha^{-n} y+r_1 w)e_2+r_1 \nu,\\
r\zeta_n &=(-\alpha^n \ell x+r_2 z)e_1+(-\alpha^{-n} \ell y+r_2 w)e_2+r_2 \nu.
\end{split}
\end{equation}
Since
\begin{equation}
\begin{split}
r^2(\eta_n^2)= & 2(\alpha^n x+r_1 z)(\alpha^{-n}y+r_1 w)+r_1^2(\nu^2),\\
r^2(\zeta_n^2)= & 2(-\alpha^n \ell x+r_2 z)(-\alpha^{-n} \ell y+r_2 w)+r_2^2 (\nu^2),\\
\end{split}
\end{equation}
We see that $(\eta_n^2)(\zeta_n^2)<0$ for $n \gg 0$ or $n \ll 0$ according as $w \ne 0$ or $z \ne 0$.

Let $W$ be a wall defined by $\w_n$.
Thus $W=(r \eta_n-r_1 \xi)^\perp$.
We take chambers ${\cal C}_\pm$ separated by $W$.
For a small neighborhood $U$ of $H \in W$, we take $H_\pm \in {\cal C}_\pm$.
Since $(r \eta_n-r_1 \xi) \cdot H=0$,
$\xi \cdot H=r(\eta_n \cdot H)/r_1=r(\zeta_n \cdot H)/r_2 \ne 0$ by the Hodge index theorem.
We may assume that 
\begin{equation}
0 \to E_1 \to E_+ \to E_2 \to 0
\end{equation}
for a general $E_+ \in { M}_{H_+}(\v)$ and
 \begin{equation}
0 \to E_2 \to E_- \to E_1 \to 0
\end{equation}
for a general $E_- \in { M}_{H_-}(\v)$,
where $E_1 \in { M}_H(\ell \w_n)$ and $E_2 \in { M}_H(\v-\ell \w_n)$. 

If $\xi \cdot H>0$, then
$h^0(E_+)h^1(E_+) \ne 0$ or
$h^0(E_-)h^1(E_-) \ne 0$ according as
$(\eta_n^2)>0$ or $(\eta_n^2)<0$.
 
If $\xi \cdot H<0$, then
$h^1(E_+)h^2(E_+) \ne 0$ or
$h^1(E_-)h^2(E_-) \ne 0$ according as
$(\eta_n^2)<0$ or $(\eta_n^2)>0$.

\end{proof}

\begin{Rem}
Lemma \ref{lem:infinite}(i) and the first claim of (ii) are well-known \cite{MW97} (see also \cite{Yamada}).
\end{Rem}

For $D \in \NS(X)$, we set 
$$
L_D:={\Bbb Z}(\xi+rD)+{\Bbb Z}(r \xi_1-r_1 \xi)\subset \NS(X).
$$
By Lemma \ref{lem:infinite},
Proposition \ref{prop:infinite} follows from the following lemma.
\begin{Lem}\label{lem:L_D}
If (i) $X$ does not contain an elliptic curve and
$\rho(X)=2$ or (ii) $\rho(X) \geq 3$, then there is a divisor $D$ such that $L_D$ is hyperbolic,
$\sqrt{-d(L_D)} \not \in {\Bbb Q}$ and
$(\xi+rD) \cdot \xi \ne 0$.
\end{Lem}

\begin{proof}
We take a divisor $D_1$ such that $L_{D_1}$ is hyperbolic and
$(\xi+r D_1) \cdot \xi \ne 0$. 
If $\sqrt{-d(L_{D_1})} \not \in {\Bbb Q}$, then the claim holds.
Hence, we  assume that $\sqrt{-d(L_{D_1})}  \in {\Bbb Q}$. Thus
$L_{D_1}$ contains isotropic elements $f$ and $g$ with $f \cdot g>0$.
We shall find a divisor $D_2$ such that
$\sqrt{-d(L_{D_1+D_2})} \not \in {\Bbb Q}$ and
$(\xi+r(D_1+D_2)) \cdot \xi \ne 0$.

We set $\xi':=\xi+rD_1$.
Then
\begin{equation}
\xi'=\frac{1}{f \cdot g}(\alpha f+\beta g),
D_1=  \frac{1}{f \cdot g}(pf+qg+h_1),
\end{equation}
where 
\begin{equation}
\alpha:=\xi' \cdot g \in {\Bbb Z},
\beta:=\xi' \cdot f \in {\Bbb Z}, 
p:=D_1 \cdot g \in {\Bbb Z},
q:=D_1 \cdot f \in {\Bbb Z},
h_1 \in L_{D_1}^\perp.
\end{equation}
We also have
$$
r_1 \xi-r \xi_1=\frac{1}{f \cdot g}(x f+y g), x,y \in {\Bbb Z}.
$$
Since $\xi' \cdot \xi \ne 0$,
$$
\xi=\frac{1}{f \cdot g}((\alpha-rp) f+(\beta-rq) g-r h_1)
$$
and $\alpha-rp \ne 0$ or $\beta-rq \ne 0$.

We set $D_2:=af+bg+h_2$, where
$a,b \in {\Bbb Z}$, $a,b \gg 0$ and $h_2 \in L_{D_1}^\perp$.
We note that $xy<0$ and $(h_2^2)<0$, and hence 
$xy(h_2^2)>0$.
\begin{equation}
\begin{split}
d(L_{D_1+D_2})=& 
\det 
\begin{pmatrix}
(\xi'+rD_2)^2 & (\xi'+rD_2) \cdot (r_1 \xi-r \xi_1)\\
 (\xi'+rD_2) \cdot (r_1 \xi-r \xi_1) &  (r_1 \xi-r \xi_1)^2
\end{pmatrix}\\
=&  (\xi'+rD_2)^2 (r_1 \xi-r \xi_1)^2- ((\xi'+rD_2) \cdot (r_1 \xi-r \xi_1))^2\\
=& -\frac{1}{(f \cdot g)^4}\left(\left((\alpha+r(f \cdot g)a)y-(\beta+r(f \cdot g)b)x \right)^2-
2r^2(h_2^2)(f \cdot g)^2xy \right).
\end{split}
\end{equation}
From our assumption on $a,b$ we get that  
$$
|((\alpha+r(f \cdot g)a)y-(\beta+r(f \cdot g)b)x|>\frac{2r^2(h_2^2)(f \cdot g)^2xy+1}{2}.
$$
Then $\sqrt{-d(L_{D_1+D_2})} \not \in {\Bbb Q}$.
Indeed, observe that if $t$ is an integer with $t>\frac{d+1}{2}$, then 
$$(t-1)^2<t^2-d<t^2.$$
Hence $\sqrt{t^2-d} \not \in {\Bbb Q}$. Letting $t=|\left((\alpha+r(f \cdot g)a)y-(\beta+r(f \cdot g)b)x \right)|$ and $d=2r^2(h_2^2)(f \cdot g)^2xy$, then gives $\sqrt{-d(L_{D_1+D_2})} \not \in {\Bbb Q}$.
Moreover, we can take $a,b$ satisfying
$$
(\xi+r(D_1+D_2))\cdot \xi=\frac{1}{(f \cdot g)^2}
\left((\beta-rq)((\alpha+r(f \cdot g)a)+(\alpha-rp)(\beta+r(f \cdot g)b) \right)
-\frac{r^2(h_1 \cdot h_2)}{f \cdot g}
\ne 0.
$$
Therefore the claim holds.
\end{proof}

\section{Ulrich bundles on generic abelian varieties}\label{sec-Ulrich}
In this section, we use weak Brill-Noether to study Ulrich bundles on abelian surfaces. 
\begin{Def}
A coherent sheaf $E$ on a projective variety $X\subset\P$ is called \emph{Ulrich} (with respect to $H:=c_1\prs{\OO_X(1)}$) if $H^\bullet\prs{X,E(-pH)}=0$ for $1\leq p\leq\dim X$. 
\end{Def}

\begin{Thm}\label{thm-ulrich}
Let $X$ be an abelian surface with an ample polarization $H$. Let $E$ be an Ulrich bundle of rank $r\geq 2$. Then its Mukai $\v$ is of the form $\v= (r, \xi, 0)e^H$ with
$$  2\xi \cdot H = r H^2, \quad \xi^2 \geq 0.$$
Conversely, suppose that $\v$ is a Mukai vector satisfying these properties and $\v'=(r, \xi, 0)$ does not admit an isotropic decomposition as in Theorem \ref{Thm:MainTheorem1}, then the general bundle $E \in M_{X, H} (\v)$ is an Ulrich bundle for $(X, H)$.
\end{Thm}

\begin{proof}
    If $E$ is an Ulrich bundle, then $h^i(E(-H)) = h^i(E(-2H))=0$ for all $i$. In particular, $\chi(E(-H))= \chi(E(-2H))=0$. Let  $\v(E(-H)) = (r, \xi, a)$. Then $\chi(E(-H))=0$ is equivalent to $a=0$. Similarly, $$\v(E(-2H))= \left(r, \xi - rH, r \frac{H^2}{2} - \xi \cdot H\right).$$ Hence $\chi(E(-2H)))=0$ if and only if 
    $2 \xi \cdot H = r H^2$. Moreover, since $\v^2 \geq 0$, we conclude that $\xi^2 \geq 0$.

    Conversely, suppose $\v = (r, \xi, 0)e^H$ satisfying $2 \xi \cdot H = r H^2$ and $\xi^2\geq 0$. Then the moduli spaces $M_{X, H}(\v')$ and $M_{X, H}(\v' e^{-H})$ are nonempty. Observe that since $\v'\ne (r, 0, -1) e^{\frac{\xi}{r}}$, if $\v'$ does not admit an isotropic decomposition, then the general bundle in  $M_{X, H}(\v')$ satisfies weak Brill-Noether and has no cohomology. Similarly for $M_{X, H}(\v' e^{-H})$. Hence, we conclude that the general bundle in $M_{X, H}(\v)$ is an Ulrich bundle.
   \end{proof}

In particular, we obtain the following corollary originally due to Beauville \cite{Be}.

\begin{Cor}\label{cor-ulrich2}
    Let $X$ be an abelian surface with an ample polarization $H$ such that $H^2 \geq 4$. Then there exists an Ulrich bundle of every even rank on $X$. Furthermore, if $H$ is divisible by $2$  in the Picard group, then there exists an Ulrich bundle of every rank $r>1$.
\end{Cor}

\begin{proof}
    Take $\v = \left(r, \frac{r}{2}H, 0\right) e^H$ in Theorem \ref{thm-ulrich}. If $r$ is even and at least 4, or if $H$ is divisible by 2 in $\Pic(X)$, then $\v$ is not primitive. Hence, by \cref{Prop:NonprimitiveCase} cannot admit an isotropic decomposition. Therefore, by Theorem \ref{thm-ulrich}, the general bundle in $M_{X, H}(\v)$ is an Ulrich bundle. Note that the dimension of this moduli space is $\frac{r^2 H^2}{4} + 2>0$.

    If $r=2$, then $\v' = (2, H, 0)$. If this vector admits an isotropic decomposition, then $\v_1 = (1, \xi_1, a)$ and $\v_2 = (1, H - \xi_1, -a)$. Using the relations $\v_1^2=\v_2^2=0$ and $\langle \v_1, \v_2 \rangle$, we calculate that $H^2 =2$ contrary to the assumption that $H^2 \geq 4$. Hence, by Theorem \ref{thm-ulrich}, the general bundle in $M_{X, H}(\v)$ is an Ulrich bundle. Note that the dimension of this moduli space is $H^2 + 2 \geq 6$.
\end{proof}

\begin{Rem}
An abelian surface cannot be embedded in $\P^n$ for $n < 4$. Hence, for a very ample $H$ on $X$, we have $\chi(\OO_X(H))= h^0(\OO_X(H))= \frac{H^2}{2}$.  We conclude that $H^2 \geq 10$, so Corollary \ref{cor-ulrich2} applies.
\end{Rem}

\begin{Cor}
    Let $X$ be an abelian surface.
    \begin{enumerate}
        \item If $\rho(X) = 1$, then there exists an Ulrich bundle of rank $r$ with respect to $mH$ if and only if $2\mid rm$.  Moreover, when an Ulrich bundle of rank $r$ exists, it has Mukai vector $\v=\left(r,\left(\frac{3rm}{2}\right)H,rm^2H^2\right)$.
        \item If $\rho(X)=2$ and $X$ contains an elliptic curve, then there is a Ulrich bundle $E$ with $\v(E)=\v$ if and only if 
$\v=(r,\xi,0)e^H$, $r \geq 2$, $\xi^2 \geq 0$ and $2 \xi \cdot H=r H^2$.
    \end{enumerate}
 \end{Cor}

 \begin{proof}
     When $\rho(X)=1$ or $\rho(X)=2$ and $X$ contains an elliptic curve, by Theorem \ref{Thm:MainTheorem1}, $\v'$ does not admit an isotropic decomposition. Hence, the corollary follows from Theorem \ref{thm-ulrich}.
 \end{proof}

\begin{Ex}
A Mukai vector $\v'= (r, \xi, 0)$ satisfying $\xi^2>0$ and $2\xi \cdot H = rH^2$ with respect to an ample $H$ may admit an isotropic decomposition. Let $X= E \times E$, where $E$ is a general elliptic curve. Then $\NS(X)$ is generated by the fiber classes $E_1, E_2$ of the two  projections and the diagonal $D$. The intersection pairing is given by $$E_1^2=E_2^2=D^2 =0, \quad E_1\cdot E_2= E_1 \cdot D= E_2\cdot D=1.$$ Let $\v'= (r, \xi, 0) = (6, -5E_1+18E_2+7D, 0)$ and  $H= 2E_1 + 5E_2$. We have $$\xi^2 = 2 \quad \mbox{and} \quad 3 (2E_1+5E_2)^2 = 60 = (2E_1 + 5E_2)(-5E_1+ 18E_2+7D).$$ Set $\v_1= (3, -2E_1+ 9E_2 + 3D, 1)$ and $\v_2= (3, -3E_1+9E_2 + 4D,-1)$. Observe that $\v_1^2=\v_2^2=0$ and $\langle \v_1, \v_2 \rangle=1$. Hence, $\v$ admits the isotropic decomposition $\v= \v_1 + \v_2$. Consequently, there are no Ulrich bundles for $H$ in the moduli space $M_{X,H}(\v' e^H)$. Note that in this example $H$ is ample with $H^2=20$, but not very ample since the degree of $H$ on $E_2$ is only 2.
\end{Ex}

\bibliographystyle{plain}

\end{document}